\documentclass[11pt,reqno]{amsart}

\title{Discrete Hamilton--Jacobi Theory}
\date{\today}
\author{Tomoki Ohsawa}
\address{Department of Mathematics, University of California, San Diego, 9500 Gilman Drive, La Jolla, California 92093--0112}
\email{tohsawa@ucsd.edu}
\author{Anthony M.~Bloch}
\address{Department of Mathematics, University of Michigan, 530 Church Street, Ann Arbor, Michigan 48109--1043}
\email{abloch@umich.edu}
\author{Melvin Leok}
\address{Department of Mathematics, University of California, San Diego, 9500 Gilman Drive, La Jolla, California 92093--0112}
\email{mleok@math.ucsd.edu}



\usepackage[margin=1in]{geometry}
\usepackage{amsthm}
\usepackage{graphicx,bigstrut,multirow,paralist}
\usepackage[numbers,sort&compress]{natbib}
\bibpunct{[}{]}{;}{n}{,}{,}

\usepackage[colorlinks=false]{hyperref}
\usepackage[all]{xy}

\theoremstyle{plain}
\newtheorem{theorem}{Theorem}[section]
\newtheorem{proposition}[theorem]{Proposition}
\newtheorem{corollary}[theorem]{Corollary}

\theoremstyle{definition}
\newtheorem{definition}[theorem]{Definition}
\newtheorem{example}[theorem]{Example}
\newtheorem{problem}[theorem]{Problem}

\theoremstyle{remark}
\newtheorem{remark}[theorem]{Remark}

\numberwithin{equation}{section}

\def\FL{\mathbb{F}L}

\def\pd#1#2{\dfrac{\partial #1}{\partial #2}}
\def\tpd#1#2{\partial #1/\partial #2}

\def\parentheses#1{\!\left(#1\right)}
\def\braces#1{\!\left\{#1\right\}}
\def\brackets#1{\!\left[#1\right]}

\def\Graph{\mathop{\mathrm{graph}}\nolimits}
\def\ext{\mathop{\mathrm{ext}}}

\def\ip#1#2{\left\langle#1,#2\right\rangle}

\def\DS{\displaystyle}
\def\setdef#1#2{ \left\{ #1 \ |\ #2 \right\} }
\def\N{\mathbb{N}}

\def\R{\mathbb{R}}

\def\defeq{\mathrel{\mathop:}=}
\def\eqdef{=\mathrel{\mathop:}}

\begin{document}
\footskip=.5in

%\allowdisplaybreaks

\begin{abstract}
  We develop a discrete analogue of Hamilton--Jacobi theory in the framework of discrete Hamiltonian mechanics.
  The resulting discrete Hamilton--Jacobi equation is discrete only in time.
  We describe a discrete analogue of Jacobi's solution and also prove a discrete version of the geometric Hamilton--Jacobi theorem.
  The theory applied to discrete linear Hamiltonian systems yields the discrete Riccati equation as a special case of the discrete Hamilton--Jacobi equation.
  We also apply the theory to discrete optimal control problems, and recover some well-known results, such as the Bellman equation (discrete-time HJB equation) of dynamic programming and its relation to the costate variable in the Pontryagin maximum principle.
  This relationship between the discrete Hamilton--Jacobi equation and Bellman equation is exploited to derive a generalized form of the Bellman equation that has controls at internal stages.
\end{abstract}

\maketitle

\section{Introduction}
\subsection{Discrete Mechanics}
Discrete mechanics is a reformulation of Lagrangian and Hamiltonian mechanics with discrete time, as opposed to a discretization of the equations in the continuous-time theory.
It not only provides a systematic view of structure-preserving integrators, but also has interesting theoretical aspects analogous to continuous-time Lagrangian and Hamiltonian mechanics~\citep[see, e.g.,][]{MaWe2001, Su2003, Su2004}.
The main feature of discrete mechanics is its use of discrete versions of variational principles.
Namely, discrete mechanics assumes that the dynamics is defined at discrete times from the outset, formulates a discrete variational principle for such dynamics, and then derives a discrete analogue of the Euler--Lagrange or Hamilton's equations from it.

The advantage of this construction is that it naturally gives rise to discrete analogues of the concepts and ideas in continuous time that have the same or similar properties, such as symplectic forms, the Legendre transformation, momentum maps, and Noether's theorem~\citep{MaWe2001}.
This in turn provides us with the discrete ingredients that facilitate further theoretical developments, such as discrete analogues of the theories of complete integrability~\citep[see, e.g.,][]{MoVe1991, Su2003, Su2004} and also those of reduction and connections~\citep{JaLeMaWe2006, LeMaWe2004, Le2004}.
Whereas the main topic in discrete mechanics is the development of structure-preserving algorithms for Lagrangian and Hamiltonian systems~\citep[see, e.g.,][]{MaWe2001}, the theoretical aspects of it are interesting in their own right, and furthermore provide insight into the numerical aspects as well.

Another notable feature of discrete mechanics, especially on the Hamiltonian side, is that it is a generalization of (nonsingular) discrete optimal control problems.
In fact, as stated in \citet{MaWe2001}, discrete mechanics is inspired by discrete formulations of optimal control problems (see, e.g., \citet{JoPo1964} and \citet{Ca1970}).

\subsection{Hamilton--Jacobi Theory}
In classical mechanics~\citep[see, e.g.,][]{La1986, Ar1991, MaRa1999, GoPoSa2001}, the Hamilton--Jacobi equation is first introduced as a partial differential equation that the action integral satisfies.
Specifically, let $Q$ be a configuration space and $T^{*}Q$ be its cotangent bundle; and let $q \in Q$ and $t > 0$ be arbitrary and suppose that $(\hat{q}(s), \hat{p}(s)) \in T^{*}Q$ is a solution of Hamilton's equations
\begin{equation}
  \label{eq:HamiltonsEq}
  \dot{q} = \pd{H}{p},
  \qquad
  \dot{p} = -\pd{H}{q}
\end{equation}
with the endpoint condition $\hat{q}(t) = q$.
Then calculate the action integral along the solution over the time interval $[0, t]$, i.e.,
\begin{equation}
  \label{eq:S(q,t)}
  S(q, t) \defeq \int_{0}^{t} \brackets{ \hat{p}(s) \cdot \dot{\hat{q}}(s) - H(\hat{q}(s), \hat{p}(s)) } ds,
\end{equation}
where we regard the resulting integral as a function of the endpoint $(q, t) \in Q \times \R_{+}$ with $\R_{+}$ being the set of positive real numbers.
Then by taking variation of the endpoint $(q, t)$, one obtains a partial differential equation satisfied by $S(q,t)$:
\begin{equation}
  \label{eq:HJEq}
  \pd{S}{t} + H\parentheses{ q, \pd{S}{q} } = 0.
\end{equation}
This is the {\em Hamilton--Jacobi equation}.

Conversely, it is shown that if $S(q,t)$ is a solution of the Hamilton--Jacobi equation then $S(q,t)$ is a generating function for the family of canonical transformations (or symplectic flows) that describe the dynamics defined by Hamilton's equations.
This result is the theoretical basis for the powerful technique of exact integration called separation of variables.

\subsection{Connection with Optimal Control and The Hamilton--Jacobi--Bellman Equation}
\label{sec:ConnectionWithOC}
The idea of Hamilton--Jacobi theory is also useful in optimal control theory~(see, e.g., \citet{Ju1997} and \citet{Be2005}).
Consider a typical optimal control problem
\begin{equation*}
  \min_{u(\cdot)} \int_{0}^{T} C(q, u)\,dt,
\end{equation*}
subject to the constraints,
\begin{equation*}
  \dot{q} = f(q, u),
\end{equation*}
and $q(0) = q_{0}$ and $q(T) = q_{T}$.
We define the augmented cost functional:
\begin{equation*}
  \hat{S}[u] \defeq \int_{0}^{T} \braces{ C(q, u) + p [ \dot{q} - f(q, u) ] } dt
  = \int_{0}^{T} \brackets{ p \cdot \dot{q} - \hat{H}(q, p, u) } dt,
\end{equation*}
where we introduced the costate $p$, and also defined the control Hamiltonian,
\begin{equation*}
  \hat{H}(q, p, u) \defeq p \cdot f(q, u) - C(q, u).
\end{equation*}
Assuming that
\begin{equation*}
  \pd{\hat{H}}{u}(q, p, u) = 0
\end{equation*}
uniquely defines the optimal control $u = u^{*}(q,p)$, we set
\begin{equation*}
  H(q,p) \defeq \max_{u}\hat{H}(q, p, u) = \hat{H}\parentheses{ q, p, u^{*}(q,p) }.
\end{equation*}
We also define the optimal cost-to-go function
\begin{align*}
  J(q, t) &\defeq \int_{t}^{T} \braces{ C(\hat{q}, u^{*}) + p \brackets{ \dot{\hat{q}} - f(\hat{q}, u^{*}) } } ds
  \\
  &= \int_{t}^{T} \brackets{ \hat{p} \cdot \dot{\hat{q}} - H(\hat{q}, \hat{p}) } ds
  = S^{*} - S(q,t),
\end{align*}
where $(\hat{q}(s), \hat{p}(s))$ for $s \in [0, T]$ is the solution of Hamilton's equations with the above $H$ such that $\hat{q}(t) = q$; and $S^{*}$ is the optimal cost
\begin{equation*}
  S^{*} \defeq \int_{0}^{T} \brackets{ \hat{p} \cdot \dot{\hat{q}} - H(\hat{q}, \hat{p}) } ds
  = \int_{0}^{T} \brackets{ \hat{p} \cdot \dot{\hat{q}} - \hat{H}\bigl( \hat{q}, \hat{p}, u^{*}(\hat{q}, \hat{p}) \bigr) } ds
  = \hat{S}[u^{*}],
\end{equation*}
and the function $S(q,t)$ is defined by
\begin{equation*}
  S(q,t) \defeq \int_{0}^{t} \brackets{ \hat{p} \cdot \dot{\hat{q}} - H(\hat{q}, \hat{p}) } ds.
\end{equation*}
Since this definition coincides with Eq.~\eqref{eq:S(q,t)}, the function $S(q,t) = S^{*} - J(q,t)$ satisfies the H--J equation~\eqref{eq:HJEq}; this reduces to the Hamilton--Jacobi--Bellman (HJB) equation for the optimal cost-to-go function $J(q,t)$:
\begin{equation}
  \label{eq:HJBEq}
  \pd{J}{t} + \min_{u}\brackets{ \pd{J}{q} \cdot f(q,u) + C(q,u) } = 0.
\end{equation}
It can also be shown that the costate $p$ of the optimal solution is related to the solution of the HJB equation.

\subsection{Discrete Hamilton--Jacobi Theory}
The main objective of this paper is to present a discrete analogue of Hamilton--Jacobi theory within the framework of discrete Hamiltonian mechanics~\cite{LaWe2006}, and also to apply the theory to discrete optimal control problems.

There are some previous works on discrete-time analogues of the Hamilton--Jacobi equation, such as \citet{ElSc1996} and \citet{LaWe2006}.
Specifically, \citet{ElSc1996} derived an equation for a generating function of a coordinate transformation that trivializes the dynamics.
This derivation is a discrete analogue of the conventional derivation of the continuous-time Hamilton--Jacobi equation~\citep[see, e.g.,][Chapter~VIII]{La1986}.
\citet{LaWe2006} formulated a discrete Lagrangian analogue of the Hamilton--Jacobi equation as a separable optimization problem.

\subsection{Main Results}
Our work was inspired by the result of \citet{ElSc1996} and starts from a reinterpretation of their result in the language of discrete mechanics.
This paper further extends the result by developing discrete analogues of results in (continuous-time) Hamilton--Jacobi theory.
Namely, we formulate a discrete analogue of Jacobi's solution, which relates the discrete action sum to a solution of the discrete Hamilton--Jacobi equation.
This also provides a very simple derivation of the discrete Hamilton--Jacobi equation and exhibits a natural correspondence with the continuous-time theory.
Another important result in this paper is a discrete analogue of the Hamilton--Jacobi theorem, which relates the solution of the discrete Hamilton--Jacobi equation with the solution of the discrete Hamilton's equations.

\begin{figure}[htbp]
  \centering
  \includegraphics[width=.9\linewidth]{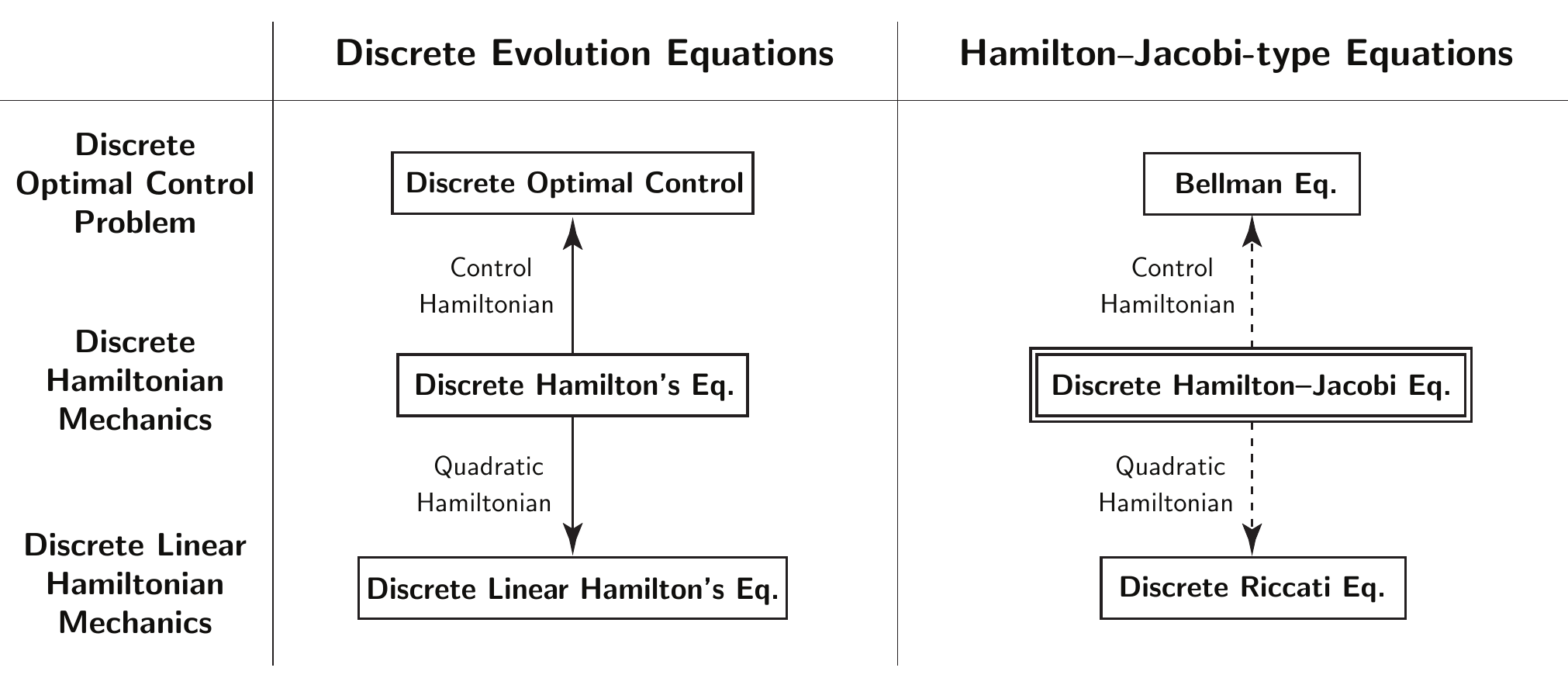}
  \caption{Discrete evolution equations (left) and corresponding discrete Hamilton--Jacobi-type equations (right).
  Dashed lines are the links established in the paper.}
  \label{fig:TextChart}
\end{figure}

We also show that the discrete Hamilton--Jacobi equation is a generalization of the discrete Riccati equation and the Bellman equation~(see Fig.~\ref{fig:TextChart}).
Specifically, we show that the discrete Hamilton--Jacobi equation applied to linear discrete Hamiltonian systems and discrete optimal control problems reduces to the discrete Riccati and Bellman equations, respectively.
This is again a discrete analogue of the well-known results that the Hamilton--Jacobi equation applied to linear Hamiltonian systems and optimal control problems reduces to the Riccati~(see, e.g., \citet[p.~421]{Ju1997}) and HJB equations~(see Section~\ref{sec:ConnectionWithOC} above), respectively.

The link between the discrete Hamilton--Jacobi equation and the Bellman equation turns out to be useful in deriving a class of generalized Bellman equations that are higher-order approximations of the original continuous-time problem.
Specifically, we use the idea of the Galerkin Hamiltonian variational integrator of \citet{LeZh2011} to derive discrete control Hamiltonians that yield higher-order approximations, and then show that the corresponding discrete Hamilton--Jacobi equation gives a class of Bellman equations with controls at internal stages.

\subsection{Outline of the Paper}
We first present a brief review of discrete Lagrangian and Hamiltonian mechanics in Section~\ref{sec:DiscreteMechanics}.
In Section~\ref{sec:DiscreteHJEq}, we describe a reinterpretation of the result of \citet{ElSc1996} in the language of discrete mechanics and a discrete analogue of Jacobi's solution to the discrete Hamilton--Jacobi equation.
The remainder of Section~\ref{sec:DiscreteHJEq} is devoted to more detailed studies of the discrete Hamilton--Jacobi equation: its left and right variants, more explicit forms of them, and also a digression on the Lagrangian side.
In Section~\ref{sec:DiscreteHJThm}, we prove a discrete version of the Hamilton--Jacobi theorem.
In Section~\ref{sec:LinearHamiltonianSystems}, we apply the theory to linear discrete Hamiltonian systems, and show that the discrete Riccati equation follows from the discrete Hamilton--Jacobi equation.
Section~\ref{sec:RelationToDTHJBEq} establishes the link with discrete-time optimal control and interprets the results of the preceding sections in this setting.
Section~\ref{sec:Higher-Order_Bellman} further extends this idea to derive a class of Bellman equations with controls at internal stages.

%We then take a harmonic oscillator as a simple physical example, and solve the discrete Hamilton--Jacobi equation explicitly.
%Finally, Section~\ref{sec:ContinuousLimit} discusses the continuous-time limit of the theory.

\section{Discrete Mechanics}
\label{sec:DiscreteMechanics}
This section briefly reviews some key results of discrete mechanics following \citet{MaWe2001} and \citet{LaWe2006}.
\subsection{Discrete Lagrangian Mechanics}
A discrete Lagrangian flow $\{ q_{k} \}_{k=0}^N$, on an $n$-dimensional differentiable manifold $Q$, can be described by the following discrete variational principle:
Let $S_{\rm d}^{N}$ be the following action sum of the discrete Lagrangian $L_{\rm d}: Q \times Q \to \R$:
\begin{equation}
  \label{eq:S_d}
  S_{\rm d}^{N}( \{q_{k}\}_{k=0}^{N} ) \defeq \sum_{k=0}^{N-1} L_{\rm d}(q_{k}, q_{k+1}) \approx \int_{0}^{t_{N}} L(q(t),\dot{q}(t))\,dt,
\end{equation}
which is an approximation of the action integral as shown above.

Consider discrete variations $q_{k} \mapsto q_{k} + \varepsilon\,\delta q_{k}$, for $k = 0,1, \dots, N$, with $\delta q_{0} = \delta q_{N} = 0$.
Then, the discrete variational principle $\delta S_{\rm d}^{N} = 0$ gives the discrete Euler--Lagrange equations:
\begin{equation}
  \label{eq:DEL}
  D_{2}L_{\rm d}(q_{k-1}, q_{k}) + D_{1}L_{\rm d}(q_{k}, q_{k+1}) = 0. 
\end{equation}
This determines the discrete flow $F_{L_{\rm d}}: Q \times Q \to Q \times Q$:
\begin{equation}
  F_{L_{\rm d}}: (q_{k-1}, q_{k}) \mapsto (q_{k}, q_{k+1}).
\end{equation}
Let us define the discrete Lagrangian symplectic one-forms $\Theta_{L_{\rm d}}^{\pm}: Q \times Q \to T^{*}(Q \times Q)$ by
\begin{subequations}
  \label{eq:DLSOF}
  \begin{align}
    \Theta_{L_{\rm d}}^{+}&: (q_{k}, q_{k+1}) \mapsto D_{2}L_{\rm d}(q_{k}, q_{k+1})\,dq_{k+1},
    \\
    \Theta_{L_{\rm d}}^{-}&: (q_{k}, q_{k+1}) \mapsto -D_{1}L_{\rm d}(q_{k}, q_{k+1})\,dq_{k}.
  \end{align}
\end{subequations}
Then, the discrete flow $F_{L_{\rm d}}$ preserves the discrete Lagrangian symplectic form
\begin{equation}
  \label{eq:DLSTF}
  \Omega_{L_{\rm d}}(q_{k}, q_{k+1}) \defeq d\Theta_{L_{\rm d}}^{+} = d\Theta_{L_{\rm d}}^{-} = D_{1}D_{2}L_{\rm d}(q_{k}, q_{k+1})\,dq_{k} \wedge dq_{k+1}.
\end{equation}
Specifically, we have
\begin{equation*}
  (F_{L_{\rm d}})^{*} \Omega_{L_{\rm d}} = \Omega_{L_{\rm d}}.
\end{equation*}

\subsection{Discrete Hamiltonian Mechanics}
Introduce the {\em right and left discrete Legendre transforms} $\FL_{\rm d}^{\pm}: Q \times Q \to T^{*}Q$ by
\begin{subequations}
  \label{eq:DLT}
  \begin{align}
    \label{eq:DLT+}
    \FL_{\rm d}^{+}&: (q_{k}, q_{k+1}) \mapsto (q_{k+1}, D_{2}L_{\rm d}(q_{k}, q_{k+1}) ),
    \\
    \label{eq:DLT-}
    \FL_{\rm d}^{-}&: (q_{k}, q_{k+1}) \mapsto (q_{k}, -D_{1}L_{\rm d}(q_{k}, q_{k+1}) ),
  \end{align}
\end{subequations}
respectively.
Then we find that the discrete Lagrangian symplectic forms Eq.~\eqref{eq:DLSOF} and \eqref{eq:DLSTF} are pull-backs by these maps of the standard symplectic form on $T^{*}Q$:
\begin{equation*}
  \Theta_{L_{\rm d}}^{\pm} = (\FL_{\rm d}^{\pm})^{*} \Theta,
  \qquad
  \Omega_{L_{\rm d}}^{\pm} = (\FL_{\rm d}^{\pm})^{*} \Omega.
\end{equation*}
Let us define the momenta
\begin{equation*}
  p_{k,k+1}^{-} \defeq -D_{1}L_{\rm d}(q_{k}, q_{k+1}),
  \qquad
  p_{k,k+1}^{+} \defeq D_{2}L_{\rm d}(q_{k}, q_{k+1}).
\end{equation*}
Then, the discrete Euler--Lagrange equations \eqref{eq:DEL} become simply $p_{k-1,k}^{+} = p_{k,k+1}^{-}$.
So defining 
\begin{equation*}
  p_{k} \defeq p_{k-1,k}^{+} = p_{k,k+1}^{-},
\end{equation*}
one can rewrite the discrete Euler--Lagrange equations \eqref{eq:DEL} as follows:
\begin{equation}
  \label{eq:DEL-H}
  \begin{array}{l}
    \DS p_{k} = -D_{1}L_{\rm d}(q_{k}, q_{k+1}),
    \medskip\\
    \DS p_{k+1} = D_{2}L_{\rm d}(q_{k}, q_{k+1}).
  \end{array}
\end{equation}
Furthermore, define the {\em discrete Hamiltonian map} $\tilde{F}_{L_{\rm d}}: T^{*}Q \to T^{*}Q$ by 
\begin{equation}
  \label{eq:DHM}
  \tilde{F}_{L_{\rm d}}: (q_{k}, p_{k}) \mapsto (q_{k+1}, p_{k+1}).
\end{equation}
Then, one may relate this map with the discrete Legendre transforms in Eq.~\eqref{eq:DLT} as follows:
\begin{equation}
  \label{eq:DHM-DLT}
  \tilde{F}_{L_{\rm d}} = \FL_{\rm d}^{+} \circ (\FL_{\rm d}^{-})^{-1}.
\end{equation}
Furthermore, one can also show that this map is symplectic, i.e.,
\begin{equation*}
  (\tilde{F}_{L_{\rm d}})^{*} \Omega = \Omega.
\end{equation*}
This is the Hamiltonian description of the dynamics defined by the discrete Euler--Lagrange equation~\eqref{eq:DEL} introduced by \citet{MaWe2001}.
Notice, however, that no discrete analogue of Hamilton's equations is introduced here, although the flow is now on the cotangent bundle $T^{*}Q$.

\citet{LaWe2006} pushed this idea further to give discrete analogues of Hamilton's equations:
From the point of view that a discrete Lagrangian is essentially a generating function of type one~\citep[][]{GoPoSa2001}, we can apply Legendre transforms to the discrete Lagrangian to find the corresponding generating functions of type two or three~\citep[][]{GoPoSa2001}.
In fact, they turn out to be a natural Hamiltonian counterpart to the discrete Lagrangian mechanics described above.
Specifically, with the right discrete Legendre transform
\begin{equation}
  p_{k+1} = \FL_{\rm d}^{+}(q_{k}, q_{k+1}) = D_{2}L_{\rm d}(q_{k}, q_{k+1}),
\end{equation}
we can define the following {\em right discrete Hamiltonian}:
\begin{equation}
  \label{eq:RDH}
  H_{\rm d}^{+}(q_{k}, p_{k+1}) = p_{k+1} \cdot q_{k+1} - L_{\rm d}(q_{k}, q_{k+1}).
\end{equation}
Then, the discrete Hamiltonian map $\tilde{F}_{L_{\rm d}}: (q_{k}, p_{k}) \mapsto (q_{k+1}, p_{k+1})$ is defined implicitly by the {\em right discrete Hamilton's equations}
\begin{subequations}
  \label{eq:RDHEq}
  \begin{align}
      \DS q_{k+1} = D_{2}H_{\rm d}^{+}(q_{k}, p_{k+1}),
      \label{eq:RDHEq-1}
      \\
      \DS p_{k} = D_{1}H_{\rm d}^{+}(q_{k}, p_{k+1}),
      \label{eq:RDHEq-2}
  \end{align}
\end{subequations}
which are precisely the characterization of a symplectic map in terms of a generating function, $H_{\rm d}^+$, of type two.
Similarly, with the left discrete Legendre transform
\begin{equation}
  p_{k} = \FL_{\rm d}^{-}(q_{k}, q_{k+1}) = -D_{1}L_{\rm d}(q_{k}, q_{k+1}),
\end{equation}
we can define the following {\em left discrete Hamiltonian}:
\begin{equation}
  \label{eq:LDH}
  H_{\rm d}^{-}(p_{k}, q_{k+1}) = -p_{k} \cdot q_{k} - L_{\rm d}(q_{k}, q_{k+1}).
\end{equation}
Then, we have the {\em left discrete Hamilton's equations}
\begin{subequations}
  \label{eq:LDHEq}
  \begin{align}
    \label{eq:LDHEq-1}
    \DS q_{k} = -D_{1}H_{\rm d}^{-}(p_{k}, q_{k+1}),
    \\
    \label{eq:LDHEq-2}
    \DS p_{k+1} = -D_{2}H_{\rm d}^{-}(p_{k}, q_{k+1}),
  \end{align}
\end{subequations}
which corresponds to a symplectic map expressed in terms of a generation function, $H_{\rm d}^-$, of type three.

On the other hand, \citet{LeZh2011} demonstrate that discrete Hamiltonian mechanics can be obtained as a direct variational discretization of continuous Hamiltonian mechanics, instead of having to go via discrete Lagrangian mechanics.

\section{Discrete Hamilton--Jacobi Equation}
\label{sec:DiscreteHJEq}
\subsection{Derivation by Elnatanov and Schiff}
\citet{ElSc1996} derived a discrete Hamilton--Jacobi equation based on the idea that the Hamilton--Jacobi equation is an equation for a symplectic change of coordinates under which the dynamics becomes trivial.
In this section, we would like to reinterpret their derivation in the framework of discrete Hamiltonian mechanics reviewed above.

\begin{theorem}
  Suppose that the discrete dynamics $\{(q_{k}, p_{k})\}_{k=0}^{N}$ is governed by the right discrete Hamilton's equations \eqref{eq:RDHEq}.
  Consider the symplectic coordinate transformation $(q_{k}, p_{k}) \mapsto (\hat{q}_{k}, \hat{p}_{k})$ that satisfies the following:
  \begin{enumerate}
    \renewcommand{\theenumi}{\roman{enumi}}
    \renewcommand{\labelenumi}{(\theenumi)}
  \item The old and new coordinates are related by the type-one generating function\footnote{This is essentially the same as Eq.~\eqref{eq:DEL-H} in the sense that they are both transformations defined by generating functions of type one: Replace $(q_{k},p_{k}, q_{k+1},p_{k+1}, L_{\rm d})$ by $(\hat{q}_{k},\hat{p}_{k}, q_{k},p_{k}, S^{k})$. However they have different interpretations: Eq.~\eqref{eq:DEL-H} describes the dynamics or time evolution whereas Eq.~\eqref{eq:COCbyS} is a change of coordinates.} $S^{k}: \R^{n} \times \R^{n} \to \R$:
    \begin{equation}
      \label{eq:COCbyS}
      \begin{array}{l}
        \DS \hat{p}_{k} = -D_{1}S^{k}(\hat{q}_{k}, q_{k}),
      \medskip\\
      \DS p_{k} = D_{2}S^{k}(\hat{q}_{k}, q_{k});
    \end{array}
  \end{equation}
  \item the dynamics in the new coordinates $\{(\hat{q}_{k}, \hat{p}_{k})\}_{k=0}^{N}$ is rendered trivial, i.e., $(\hat{q}_{k+1}, \hat{p}_{k+1}) = (\hat{q}_{k}, \hat{p}_{k})$.
  \end{enumerate}
  Then, the set of functions $\{ S^{k} \}_{k=1}^{N}$ satisfies the {\em discrete Hamilton--Jacobi equation}:
  \begin{equation}
    \label{eq:RDHJEq-ElSc}
    S^{k+1}(\hat{q}_{0}, q_{k+1}) - S^{k}(\hat{q}_{0}, q_{k})
    - D_{2}S^{k+1}(\hat{q}_{0}, q_{k+1}) \cdot q_{k+1}
    + H_{\rm d}^{+}\parentheses{ q_{k}, D_{2}S^{k+1}(\hat{q}_{0}, q_{k+1}) }
    = 0,
  \end{equation}
  or, with the shorthand notation $S_{\rm d}^{k}(q_{k}) \defeq S^{k}(\hat{q}_{0}, q_{k})$, 
  \begin{equation}
    \label{eq:RDHJEq-implicit}
    S_{\rm d}^{k+1}(q_{k+1}) - S_{\rm d}^{k}(q_{k})
    - DS_{\rm d}^{k+1}(q_{k+1}) \cdot q_{k+1}
    + H_{\rm d}^{+}\parentheses{ q_{k}, DS_{\rm d}^{k+1}(q_{k+1}) }
    = 0.
  \end{equation}
\end{theorem}

\begin{proof}
  The key ingredient in the proof is the right discrete Hamiltonian in the new coordinates, i.e., a function $\hat{H}_{\rm d}^{+}(\hat{q}_{k}, \hat{p}_{k+1})$ that satisfies
  \begin{equation}
    \label{eq:RDHEq-ElSc-2a}
    \begin{array}{l}
      \DS \hat{q}_{k+1} = D_{2}\hat{H}_{\rm d}^{+}(\hat{q}_{k}, \hat{p}_{k+1}),
      \medskip\\
      \DS \hat{p}_{k} = D_{1}\hat{H}_{\rm d}^{+}(\hat{q}_{k}, \hat{p}_{k+1}),
    \end{array}
  \end{equation}
  or equivalently,
  \begin{equation}
    \label{eq:RDHEq-ElSc-2b}
    \hat{p}_{k}\,d\hat{q}_{k} + \hat{q}_{k+1}\,d\hat{p}_{k+1} = d\hat{H}_{\rm d}^{+}(\hat{q}_{k}, \hat{p}_{k+1}).
  \end{equation}

  Let us first write $\hat{H}_{\rm d}^{+}$ in terms of the original right discrete Hamiltonian $H_{\rm d}^{+}$ and the generating function $S^{k}$.
  For that purpose, first rewrite Eqs.~\eqref{eq:RDHEq} and \eqref{eq:COCbyS} as follows:
  \begin{equation*}
    p_{k}\,dq_{k} = -q_{k+1}\,dp_{k+1} + dH_{\rm d}^{+}(q_{k}, p_{k+1})
  \end{equation*}
  and
  \begin{equation*}
    \hat{p}_{k}\,d\hat{q}_{k} = p_{k}\,dq_{k} - dS^{k}(\hat{q}_{k}, q_{k}),
  \end{equation*}
  respectively.
  Then, using the above relations, we have
  \begin{align*}
    \hat{p}_{k}\,d\hat{q}_{k} + \hat{q}_{k+1}\,d\hat{p}_{k+1} &= \hat{p}_{k}\,d\hat{q}_{k} + d(\hat{p}_{k+1} \cdot \hat{q}_{k+1}) - \hat{p}_{k+1}\,d\hat{q}_{k+1}
    \\
    &= p_{k}\,dq_{k} - dS^{k}(\hat{q}_{k}, q_{k}) + d(\hat{p}_{k+1} \cdot \hat{q}_{k+1}) - p_{k+1}\,dq_{k+1} + dS^{k+1}(\hat{q}_{k+1}, q_{k+1})
    \\
    &= -q_{k+1}\,dp_{k+1} + dH_{\rm d}^{+}(q_{k}, p_{k+1})
    \\
    &\quad - dS^{k}(\hat{q}_{k}, q_{k}) + d(\hat{p}_{k+1} \hat{q}_{k+1}) - p_{k+1}\,dq_{k+1} + dS^{k+1}(\hat{q}_{k+1}, q_{k+1})
    \\
    &= d\parentheses{
      H_{\rm d}^{+}(q_{k}, p_{k+1}) + \hat{p}_{k+1} \cdot \hat{q}_{k+1} - p_{k+1} \cdot q_{k+1}
      + S^{k+1}(\hat{q}_{k+1}, q_{k+1}) - S^{k}(\hat{q}_{k}, q_{k}) 
    }.
  \end{align*}
  Thus, in view of Eq.~\eqref{eq:RDHEq-ElSc-2b}, we obtain
  \begin{equation}
    \label{eq:hatH_d+}
    \hat{H}_{\rm d}^{+}(\hat{q}_{k}, \hat{p}_{k+1})
    = H_{\rm d}^{+}(q_{k}, p_{k+1}) + \hat{p}_{k+1} \cdot \hat{q}_{k+1} - p_{k+1} \cdot q_{k+1}
    + S^{k+1}(\hat{q}_{k+1}, q_{k+1}) - S^{k}(\hat{q}_{k}, q_{k}).
  \end{equation}
  
  Now consider the choice of the new right discrete Hamiltonian $\hat{H}_{\rm d}^{+}$ that renders the dynamics trivial, i.e., $(\hat{q}_{k+1}, \hat{p}_{k+1}) = (\hat{q}_{k}, \hat{p}_{k})$.
  It is clear from Eq.~\eqref{eq:RDHEq-ElSc-2a} that we can set
  \begin{equation*}
    \hat{H}_{\rm d}^{+}(\hat{q}_{k}, \hat{p}_{k+1}) = \hat{p}_{k+1} \cdot \hat{q}_{k}.
  \end{equation*}
  Then, Eq.~\eqref{eq:hatH_d+} becomes
  \begin{equation*}
    \hat{p}_{k+1} \cdot \hat{q}_{k}
    = H_{\rm d}^{+}(q_{k}, p_{k+1}) + \hat{p}_{k+1} \cdot \hat{q}_{k+1} - p_{k+1} \cdot q_{k+1}
    + S^{k+1}(\hat{q}_{k+1}, q_{k+1}) - S^{k}(\hat{q}_{k}, q_{k}),
  \end{equation*}
  and since $\hat{q}_{k+1} = \hat{q}_{k} = \dots = \hat{q}_{0}$, we have
  \begin{equation*}
    0 = H_{\rm d}^{+}(q_{k}, p_{k+1}) - p_{k+1} \cdot q_{k+1}
    + S^{k+1}(\hat{q}_{0}, q_{k+1}) - S^{k}(\hat{q}_{0}, q_{k}).
  \end{equation*}
  Eliminating $p_{k+1}$ by using Eq.~\eqref{eq:COCbyS}, we obtain Eq.~\eqref{eq:RDHJEq-ElSc}.
\end{proof}

\begin{remark}
  What \citet{ElSc1996} refer to as the {\em Hamilton--Jacobi difference equation} is the following:
  \begin{equation}
    \label{eq:HJDEq}
    S^{k+1}(\hat{q}_{0}, q_{k+1}) - S^{k}(\hat{q}_{0}, q_{k})
    - D_{2}S^{k+1}(\hat{q}_{0}, q_{k+1}) \cdot D_{2}H_{\rm d}^{+}(q_{k}, p_{k+1})
    + H_{\rm d}^{+}(q_{k}, p_{k+1})
    = 0.
  \end{equation}
  It is clear that this is equivalent to Eq.~\eqref{eq:RDHJEq-ElSc} in view of Eq.~\eqref{eq:RDHEq}
\end{remark}

\subsection{Discrete Analogue of Jacobi's Solution}
This section presents a discrete analogue of Jacobi's solution.
This also gives an alternative derivation of the discrete Hamilton--Jacobi equation that is much simpler than the one shown above.
\begin{theorem}
  \label{thm:DJS}
  Consider the action sums, Eq.~\eqref{eq:S_d}, written in terms of the right discrete Hamiltonian, Eq.~\eqref{eq:RDH}:
  \begin{equation}
    \label{eq:S_d+}
    S_{\rm d}^{k}(q_{k})
    \defeq \sum_{l=0}^{k-1} \brackets{ p_{l+1} \cdot q_{l+1} - H_{\rm d}^{+}(q_{l}, p_{l+1}) }
  \end{equation}
  evaluated along a solution of the right discrete Hamilton's equations \eqref{eq:RDHEq}; each $S_{\rm d}^{k}(q_{k})$ is seen as a function of the end point coordinates $q_{k}$ and the discrete end time $k$.
  Then, these action sums satisfy the discrete Hamilton--Jacobi equation~\eqref{eq:RDHJEq-implicit}.
\end{theorem}

\begin{proof}
  From Eq.~\eqref{eq:S_d+}, we have
  \begin{equation}
    \label{eq:diffS_d+}
    S_{\rm d}^{k+1}(q_{k+1}) - S_{\rm d}^{k}(q_{k}) = p_{k+1} \cdot q_{k+1} - H_{\rm d}^{+}(q_{k}, p_{k+1}),
  \end{equation}
  where $p_{k+1}$ is considered to be a function of $q_{k}$ and $q_{k+1}$, i.e., $p_{k+1} = p_{k+1}(q_{k}, q_{k+1})$.
  Taking the derivative of both sides with respect to $q_{k+1}$, we have
  \begin{equation*}
    DS_{\rm d}^{k+1}(q_{k+1}) = p_{k+1} + \pd{p_{k+1}}{q_{k+1}} \cdot \brackets{ q_{k+1} - D_{2}H_{\rm d}^{+}(q_{k}, p_{k+1}) }.
  \end{equation*}
  However, the terms in the brackets vanish because the right discrete Hamilton's equations \eqref{eq:RDHEq} are assumed to be satisfied.
  Thus, we have
  \begin{equation}
    \label{eq:p-dS+}
    p_{k+1} = DS_{\rm d}^{k+1}(q_{k+1}).
  \end{equation}
  Substituting this into Eq.~\eqref{eq:diffS_d+} gives Eq.~\eqref{eq:RDHJEq-implicit}.
\end{proof}

\begin{remark}
  \label{remark:continuous_HJ}
  Recall that, in the derivation of the continuous Hamilton--Jacobi equation~\citep[see, e.g.,][Section~23]{GeFo2000}, we consider the variation of the action integral, Eq.~\eqref{eq:S(q,t)}, with respect to the end point $(q,t)$ and find
  \begin{equation}
    \label{eq:dS}
    dS = p\,dq - H(q,p)\,dt.
  \end{equation}
  This gives
  \begin{equation*}
    \pd{S}{t} = -H(q,p),
    \qquad
    p = \pd{S}{q},
  \end{equation*}
  and hence the Hamilton--Jacobi equation
  \begin{equation*}
    \pd{S}{t} + H\parentheses{ q, \pd{S}{q} } = 0.
  \end{equation*}
  In the above derivation of the discrete Hamilton--Jacobi equation \eqref{eq:RDHJEq-implicit}, the difference in two action sums, Eq.~\eqref{eq:diffS_d+}, is a natural discrete analogue of the variation $dS$ in Eq.~\eqref{eq:dS}.
  Notice also that Eq.~\eqref{eq:diffS_d+} plays the same essential role as Eq.~\eqref{eq:dS} does in deriving the Hamilton--Jacobi equation.
  Table~\ref{table:Correspondence} summarizes the correspondence between the ingredients in the continuous and discrete theories (see also Remark~\ref{remark:continuous_HJ}).
  \begin{table}[htbp]
    \centering
    \caption{Correspondence between ingredients in continuous and discrete theories; $\R_{\ge0}$ is the set of non-negative real numbers and $\N_{0}$ is the set of non-negative integers.}
    \label{table:Correspondence}
    \renewcommand{\arraystretch}{1.5}
    \small
    \begin{tabular}{|c|c|}
      \hline
      {\bf Continuous} & {\bf Discrete}
      \\\hline\hline
      $\bigstrut (q, t) \in Q \times \R_{\ge0} $ & $(q_{k}, k) \in Q \times \N_{0}$
      \\\hline
      $\DS \dot{q} = \tpd{H}{p}, $ & $\DS q_{k+1} = D_{2}H_{\rm d}^{+}(q_{k}, p_{k+1})$,
      \\
      $\DS \dot{p} = -\tpd{H}{q} $ & $\DS p_{k} = D_{1}H_{\rm d}^{+}(q_{k}, p_{k+1})$
      \\\hline
      \multirow{2}{*}{$\DS S(q, t) \defeq \int_{0}^{t} \brackets{ p(s) \cdot \dot{q}(s) - H(q(s), p(s)) } ds$}
      &
      \multirow{2}{*}{$\DS S_{\rm d}^{k}(q_{k}) \defeq \sum_{l=0}^{k-1} \brackets{ p_{l+1} \cdot q_{l+1} - H_{\rm d}^{+}(q_{l}, p_{l+1}) }$}
      \\
      &
      \\\hline
      \multirow{2}{*}{$\DS dS = \pd{\mathstrut S}{\mathstrut q}\,dq + \pd{S}{t}\,dt$}
      &
      \multirow{2}{*}{$\DS S_{\rm d}^{k+1}(q_{k+1}) - S_{\rm d}^{k}(q_{k})$}
      \\
      &
      \\\hline
      $\DS \bigstrut p\,dq - H(q,p)\,dt $ & $\DS p_{k+1} \cdot q_{k+1} - H_{\rm d}^{+}(q_{k}, p_{k+1})$
      \\\hline
      \multirow{2}{*}{$\DS \pd{S}{t} + H\parentheses{ q, \pd{S}{q} } = 0$}
      &
      $\DS S_{\rm d}^{k+1}(q_{k+1}) - S_{\rm d}^{k}(q_{k})
      - DS_{\rm d}^{k+1}(q_{k+1}) \cdot q_{k+1}$
      \\
      &
      \hfill$\DS + H_{\rm d}^{+}\parentheses{ q_{k}, DS_{\rm d}^{k+1}(q_{k+1}) } = 0$
      \\\hline
    \end{tabular}
  \end{table}
\end{remark}

\subsection{The Right and Left Discrete Hamilton--Jacobi Equations}
Recall that, in Eq.~\eqref{eq:S_d+}, we wrote the action sum, Eq.~\eqref{eq:S_d}, in terms of the right discrete Hamiltonian,~Eq.~\eqref{eq:RDH}.
We can also write it in terms of the left discrete Hamiltonian, Eq.~\eqref{eq:LDH}, as follows:
\begin{equation}
  \label{eq:S_d-}
  S_{\rm d}^{k}(q_{k})
  = \sum_{l=0}^{k-1} \brackets{ -p_{l} \cdot q_{l} - H_{\rm d}^{-}(p_{l}, q_{l+1}) }.
\end{equation}
Then, we can proceed as in the proof of Theorem~\ref{thm:DJS}:
First, we have
\begin{equation}
  \label{eq:diffS_d-}
  S_{\rm d}^{k+1}(q_{k+1}) - S_{\rm d}^{k}(q_{k}) = -p_{k} \cdot q_{k} - H_{\rm d}^{-}(p_{k}, q_{k+1}).
\end{equation}
where $p_{k}$ is considered to be a function of $q_{k}$ and $q_{k+1}$, i.e., $p_{k} = p_{k}(q_{k}, q_{k+1})$.
Taking the derivative of both sides with respect to $q_{k}$, we have
\begin{equation*}
  -DS_{\rm d}^{k}(q_{k}) = -p_{k} - \pd{p_{k}}{q_{k}} \cdot \brackets{ q_{k} + D_{1}H_{\rm d}^{-}(p_{k}, q_{k+1}) }.
\end{equation*}
However, the terms in the brackets vanish because the left discrete Hamilton's equations \eqref{eq:LDHEq} are assumed to be satisfied.
Thus, we have
\begin{equation}
  \label{eq:p-dS-}
  p_{k} = DS_{\rm d}^{k}(q_{k}).
\end{equation}
Substituting this into Eq.~\eqref{eq:diffS_d-} gives the discrete Hamilton--Jacobi equation with the left discrete Hamiltonian:
\begin{equation}
  \label{eq:LDHJEq-implicit}
  S_{\rm d}^{k+1}(q_{k+1}) - S_{\rm d}^{k}(q_{k}) 
  + DS_{\rm d}^{k}(q_{k}) \cdot q_{k}
  + H_{\rm d}^{-}\parentheses{ DS_{\rm d}^{k}(q_{k}), q_{k+1} }
  = 0.
\end{equation}
We refer to Eqs.~\eqref{eq:RDHJEq-implicit} and \eqref{eq:LDHJEq-implicit} as the {\em right and left discrete Hamilton--Jacobi equations}, respectively.

As mentioned above, Eqs.~\eqref{eq:S_d+} and \eqref{eq:S_d-} are the same action sum,~Eq.~\eqref{eq:S_d}, expressed in different ways.
Therefore we may summarize the above argument as follows:
\begin{proposition}
  The action sums, Eq.~\eqref{eq:S_d+} or equivalently Eq.~\eqref{eq:S_d-}, satisfy both the right and left discrete Hamilton--Jacobi equations,~\eqref{eq:RDHJEq-implicit} and \eqref{eq:LDHJEq-implicit}, respectively.
\end{proposition}

\subsection{Explicit Forms of the Discrete Hamilton--Jacobi Equations}
\label{ssec:ExplicitDHJEq}
The expressions for the right and left discrete Hamilton--Jacobi equations in Eqs.~\eqref{eq:RDHJEq-implicit} and \eqref{eq:LDHJEq-implicit} are implicit in the sense that they contain two spatial variables $q_{k}$ and $q_{k+1}$; Theorem~\ref{thm:DJS} suggests that one may consider $q_{k}$ and $q_{k+1}$ to be related by the discrete Hamiltonian dynamics defined by either the right or left discrete Hamilton's equations~~\eqref{eq:RDHEq} or \eqref{eq:LDHEq}, or equivalently, the discrete Hamiltonian map $\tilde{F}_{L_{\rm d}}:(q_{k}, p_{k}) \mapsto (q_{k+1}, p_{k+1})$ defined in Eq.~\eqref{eq:DHM}.
More specifically, we may write $q_{k+1}$ in terms of $q_{k}$.
This results in explicit forms of the discrete Hamilton--Jacobi equations, and we shall {\em define} the discrete Hamilton--Jacobi equations by the resulting explicit forms.
We will see later in Section~\ref{sec:RelationToDTHJBEq} that the explicit form is compatible with the formulation of the Bellman equation.

For the right discrete Hamilton--Jacobi equation~\eqref{eq:RDHJEq-implicit}, we first define the map $f^{+}_{k}: Q \to Q$ as follows:
Replace $p_{k+1}$ in Eq.~\eqref{eq:RDHEq-1} by $DS_{\rm d}^{k+1}(q_{k+1})$ as suggested by Eq.~\eqref{eq:p-dS+}:
\begin{equation*}
  q_{k+1} = D_{2}H_{\rm d}^{+}\parentheses{ q_{k}, DS_{\rm d}^{k+1}(q_{k+1}) }.
\end{equation*}
Assuming this equation is solvable for $q_{k+1}$, we define $f_{k}^{+}: Q \to Q$ by $f_k(q_k)=q_{k+1}$, i.e., $f^{+}_{k}$ is implicitly defined by
\begin{equation}
  \label{eq:f^+_k-implicit}
  f^{+}_{k}(q_{k}) = D_{2}H_{\rm d}^{+}\parentheses{ q_{k}, DS_{\rm d}^{k+1}(f^{+}_{k}(q_{k})) }.
\end{equation}
We may now identify $q_{k+1}$ with $f^{+}_{k}(q_{k})$ in the implicit form of the right Hamilton--Jacobi equation~\eqref{eq:RDHJEq-implicit}:
\begin{equation}
  \label{eq:RDHJEq}
  S_{\rm d}^{k+1}(f^{+}_{k}(q)) - S_{\rm d}^{k}(q) 
  - DS_{\rm d}^{k+1}(f^{+}_{k}(q)) \cdot f^{+}_{k}(q)
  + H_{\rm d}^{+}\parentheses{ q, DS_{\rm d}^{k+1}(f^{+}_{k}(q)) }
  = 0,
\end{equation}
where we suppressed the subscript $k$ of $q_{k}$ since it is now clear that $q_{k}$ is an independent variable as opposed to a function of the discrete time $k$.
We {\em define} Eq.~\eqref{eq:RDHJEq} to be the {\em right discrete Hamilton--Jacobi equation}.
Notice that these are differential-difference-functional equations defined on $Q \times \N$, with the spatial variable $q$ and the discrete time $k$.

For the left discrete Hamilton--Jacobi equation~\eqref{eq:LDHJEq-implicit}, we define the map $f^{-}_{k}: Q \to Q$ as follows:
\begin{equation}
  \label{eq:f^-_k}
  f^{-}_{k}(q_{k}) \defeq \pi_{Q} \circ \tilde{F}_{L_{\rm d}} \parentheses{ dS_{\rm d}^{k}(q_{k}) },
\end{equation}
where $\pi_{Q}: T^{*}Q \to Q$ is the cotangent bundle projection; equivalently, $f^{-}_{k}$ is defined so that the diagram below commutes.
\begin{equation*}
  \vcenter{
    \xymatrix@!0@R=0.75in@C=1in{
      T^{*}Q \ar[r]^{\tilde{F}_{L_{\rm d}}} & T^{*}Q \ar[d]^{\pi_{Q}}
      \\
      Q \ar[u]^{dS_{\rm d}^{k}} \ar@{-->}[r]_{f^{-}_{k}} & Q
    }
  }
  \qquad
  \vcenter{
    \xymatrix@!0@R=0.75in@C=1.1in{
      dS_{\rm d}^{k}(q_{k}) \ar@{|->}[r] & \tilde{F}_{L_{\rm d}}\parentheses{ dS_{\rm d}^{k}(q_{k}) } \ar@{|->}[d]
      \\
      q_{k} \ar@{|->}[u] \ar@{|-->}[r] & f^{-}_{k}(q_{k})
    }
  }
\end{equation*}
Notice also that, since the map $\tilde{F}_{L_{\rm d}}: (q_{k}, p_{k}) \mapsto (q_{k+1}, p_{k+1})$ is defined by Eq.~\eqref{eq:LDHEq}, $f^{-}_{k}$ is defined implicitly by
\begin{equation}
  \label{eq:f^-_k-implicit}
  q_{k} = -D_{1}H_{\rm d}^{-}\parentheses{ DS_{\rm d}^{k}(q_{k}), f^{-}_{k}(q_{k}) }.
\end{equation}
In other words, replace $p_{k}$ in Eq.~\eqref{eq:LDHEq-1} by $DS_{\rm d}^{k}(q_{k})$ as suggested by Eq.~\eqref{eq:p-dS-}, and define $f_{k}^{-}(q_{k})$ as the $q_{k+1}$ in the resulting equation.

We may now identify $q_{k+1}$ with $f^{-}_{k}(q_{k})$ in Eq.~\eqref{eq:LDHJEq-implicit}:
\begin{equation}
  \label{eq:LDHJEq}
  S_{\rm d}^{k+1}(f^{-}_{k}(q)) - S_{\rm d}^{k}(q) 
  + DS_{\rm d}^{k}(q) \cdot q
  + H_{\rm d}^{-}\parentheses{ DS_{\rm d}^{k}(q), f^{-}_{k}(q) }
  = 0,
\end{equation}
where we again suppressed the subscript $k$ of $q_{k}$.
We {\em define} Eqs.~\eqref{eq:RDHJEq} and \eqref{eq:LDHJEq} to be the {\em right and left discrete Hamilton--Jacobi equations}, respectively.

\begin{remark}
  \label{remark:complexity_of_DHJEq}
  Notice that the right discrete Hamilton--Jacobi equation~\eqref{eq:RDHJEq} is more complicated than the left one~\eqref{eq:LDHJEq}, particularly because the map $f^{+}_{k}$ appears more often than $f^{-}_{k}$ does in the latter; notice here that, as shown in Eq.~\eqref{eq:f^-_k}, the maps $f^{\pm}_{k}$ in the discrete Hamilton--Jacobi equations~\eqref{eq:RDHJEq} and \eqref{eq:LDHJEq} depend on the function $S_{\rm d}^{k}$, which is the unknown one has to solve for.

  However, it is possible to define an equally simple variant of the right discrete Hamilton--Jacobi equation by writing $q_{k-1}$ in terms of $q_{k}$: 
  Let us first define $g_{k}: Q \to Q$ by
  \begin{equation*}
    g_{k}(q_{k}) \defeq \pi_{Q} \circ \tilde{F}_{L_{\rm d}}^{-1} \parentheses{ dS_{\rm d}^{k}(q_{k}) },
  \end{equation*}
  or so that the diagram below commutes.
  \begin{equation*}
    \vcenter{
      \xymatrix@!0@R=0.75in@C=1in{
        T^{*}Q \ar[d]_{\pi_{Q}} & T^{*}Q \ar[l]_{\tilde{F}_{L_{\rm d}}^{-1}}
        \\
        Q & Q \ar[u]_{dS_{\rm d}^{k}} \ar@{-->}[l]^{g_{k}}
      }
    }
    \qquad
    \vcenter{
      \xymatrix@!0@R=0.75in@C=1.1in{
        \tilde{F}_{L_{\rm d}}^{-1}\parentheses{ dS_{\rm d}^{k}(q_{k}) } \ar@{|->}[d] & dS_{\rm d}^{k}(q_{k}) \ar@{|->}[l]
        \\
        g_{k}(q_{k}) & q_{k} \ar@{|->}[u] \ar@{|-->}[l]
      }
    }
  \end{equation*}
  Now, in Eq.~\eqref{eq:RDHJEq-implicit}, change the indices from $(k,k+1)$ to $(k-1,k)$ and identify $q_{k-1}$ with $g_{k}(q_{k})$ to obtain
  \begin{equation*}
    S_{\rm d}^{k}(q) - S_{\rm d}^{k-1}(g_{k}(q)) 
    - DS_{\rm d}^{k}(q) \cdot q
    + H_{\rm d}^{+}\parentheses{ g_{k}(q), DS_{\rm d}^{k}(q) }
    = 0,
  \end{equation*}
  where we again suppressed the subscript $k$ of $q_{k}$.
  This is as simple as the left discrete Hamilton--Jacobi equation~\eqref{eq:LDHJEq}.
  However the map $g_{k}$ is, being backward in time, rather unnatural compared to $f_{k}$.
  Furthermore, as we shall see in Section~\ref{sec:RelationToDTHJBEq}, in the discrete optimal control setting, the map $f_{k}$ is defined by a given function and thus the formulation with $f_{k}$ will turn out to be more convenient.
\end{remark}

\subsection{The Discrete Hamilton--Jacobi Equation on the Lagrangian Side}
First, notice that Eq.~\eqref{eq:S_d} gives
\begin{equation}
  \label{eq:DHJEq-L}
  S_{\rm d}^{k+1}(q_{k+1}) - S_{\rm d}^{k}(q_{k}) = L_{\rm d}(q_{k}, q_{k+1}).
\end{equation}
This is essentially the Lagrangian equivalent of the discrete Hamilton--Jacobi equation~\eqref{eq:RDHJEq} as \citet{LaWe2006} suggest.
Let us apply the same argument as above to obtain the explicit form for Eq.~\eqref{eq:DHJEq-L}.
Taking the derivative of the above equation with respect to $q_{k}$, we have
\begin{equation*}
  -D_{1}L_{\rm d}(q_{k}, q_{k+1})\,dq_{k} = dS_{\rm d}^{k}(q_{k}),
\end{equation*}
and hence from the definition of the left discrete Legendre transform,~Eq.~\eqref{eq:DLT-},
\begin{equation*}
  \FL_{\rm d}^{-}(q_{k}, q_{k+1}) = dS_{\rm d}^{k}(q_{k}).
\end{equation*}
Assuming that $\FL_{\rm d}^{-}$ is invertible, we have
\begin{equation*}
   (q_{k}, q_{k+1}) =  (\FL_{\rm d}^{-})^{-1}\parentheses{ dS_{\rm d}^{k}(q_{k}) } \eqdef (q_{k}, f^{L}_{k}(q_{k})),
\end{equation*}
where we defined the map $f^{L}_{k}: Q \to Q$ as follows (see the commutative diagram below):
\begin{equation}
  \label{eq:f_k-L}
  f^{L}_{k}(q_{k}) \defeq pr_{2} \circ (\FL_{\rm d}^{-})^{-1}\parentheses{ dS_{\rm d}^{k}(q_{k}) },
\end{equation}
where $pr_{2}: Q \times Q \to Q$ is the projection to the second factor, i.e., $pr_{2}(q_{1}, q_{2}) = q_{2}$.
Thus, eliminating $q_{k+1}$ from Eq.~\eqref{eq:DHJEq-L} and then replacing $q_{k}$ by $q$, we obtain the discrete Hamilton--Jacobi equation on the Lagrangian side:
\begin{equation}
  \label{eq:DHJEq-L2}
  S_{\rm d}^{k+1}(f^{L}_{k}(q)) - S_{\rm d}^{k}(q) = L_{\rm d}\parentheses{ q, f^{L}_{k}(q) }.
\end{equation}
The map $f^{L}_{k}$ defined in Eq.~\eqref{eq:f_k-L} is identical to $f^{-}_{k}$ defined above in Eq.~\eqref{eq:f^-_k} as the commutative diagram below demonstrates.
\begin{equation*}
  \vcenter{
    \xymatrix@!0@R=0.60in@C=.70in{
      T^{*}Q \ar[rr]^{\tilde{F}_{L_{\rm d}}} \ar[rd]^{\!\!(\FL_{\rm d}^{-})^{-1}} & & T^{*}Q \ar[dd]^{\pi_{Q}}
      \\
      & Q \times Q \ar[ru]^{\FL_{\rm d}^{+}\!\!} \ar[ld]^{pr_{1}} \ar[rd]_{pr_{2}} & 
      \\
      Q \ar[uu]^{dS_{\rm d}^{k}} \ar@{-->}[rr]_{f^{L}_{k}, f^{-}_{k}} & & Q
    }
  }
  \qquad
  \vcenter{
    \xymatrix@!0@R=0.60in@C=.70in{
      dS_{\rm d}^{k}(q_{k}) \ar@{|->}[rr] \ar@{|->}[rd] & & \tilde{F}_{L_{\rm d}}\parentheses{ dS_{\rm d}^{k}(q_{k}) } \ar@{|->}[dd]
      \\
      & (q_{k}, f^{L}_{k}(q_{k})) \ar@{|->}[ru] \ar@{|->}[ld] \ar@{|->}[rd] & 
      \\
      q_{k} \ar@{|->}[uu] \ar@{|-->}[rr] & & f^{L}_{k}(q_{k})
    }
  }
\end{equation*}
The commutativity of the square in the diagram defines the $f^{-}_{k}$ as we saw earlier, whereas that of the right-angled triangle on the lower left defines the $f^{L}_{k}$ in Eq.~\eqref{eq:f_k-L}; note the relation $\tilde{F}_{L_{\rm d}} = \FL_{\rm d}^{+} \circ (\FL_{\rm d}^{-})^{-1}$ from Eq.~\eqref{eq:DHM-DLT}.

The map $f^{L}_{k}$ being identical to $f^{-}_{k}$ implies that the discrete Hamilton--Jacobi equations on the Hamiltonian and Lagrangian sides, Eqs.~\eqref{eq:LDHJEq} and \eqref{eq:DHJEq-L2}, are equivalent.

\section{Discrete Hamilton--Jacobi Theorem}
\label{sec:DiscreteHJThm}
The following gives a discrete analogue of the geometric Hamilton--Jacobi theorem by \citet[Theorem~5.2.4]{AbMa1978}:
\begin{theorem}[Discrete Hamilton--Jacobi]
  \label{thm:DHJ}
  Suppose that $S_{\rm d}^{k}$ satisfies the right discrete Hamilton--Jacobi equation \eqref{eq:RDHJEq}, and let $\{ c_{k} \}_{k=0}^{N} \subset Q$ be a set of points such that
  \begin{equation}
    \label{eq:RDHJ-curve}
    c_{k+1} = f^{+}_{k}(c_{k})
    \quad
    \text{for}
    \quad
    k = 0, 1, \dots, N-1.
  \end{equation}
  Then, the set of points $\{ (c_{k}, p_{k}) \}_{k=0}^{N} \subset T^{*}Q$ with
  \begin{equation}
    \label{eq:RDHJ_p-dS}
    p_{k} \defeq DS_{\rm d}^{k}(c_{k})
  \end{equation}
  is a solution of the right discrete Hamilton's equations \eqref{eq:RDHEq}.
  
  Similarly, suppose that $S_{\rm d}^{k}$ satisfies the left discrete Hamilton--Jacobi equation \eqref{eq:LDHJEq}, and let $\{ c_{k} \}_{k=0}^{N} \subset Q$ be a set of points that satisfy
  \begin{equation}
    \label{eq:LDHJ-curve}
    c_{k+1} = f^{-}_{k}(c_{k})
    \quad
    \text{for}
    \quad
    k = 0, 1, \dots, N-1.
  \end{equation}
  Furthermore, assume that the Jacobian $Df^{-}_{k}$ is invertible at each point $c_{k}$.
  Then, the set of points $\{ (c_{k}, p_{k}) \}_{k=0}^{N} \subset T^{*}Q$ with 
  \begin{equation}
    \label{eq:LDHJ_p-dS}
    p_{k} \defeq DS_{\rm d}^{k}(c_{k})
  \end{equation}
  is a solution of the left discrete Hamilton's equations \eqref{eq:LDHEq}.
\end{theorem}

\begin{proof}
  To prove the first assertion, first recall the implicit definition of $f^{+}_{k}$ in Eq.~\eqref{eq:f^+_k-implicit}:
  \begin{equation}
    \label{eq:f_k-D2Hd+}
    f^{+}_{k}(q) = D_{2}H_{\rm d}^{+}\parentheses{ q, DS_{\rm d}^{k+1}(f^{+}_{k}(q)) }.
  \end{equation}
  In particular, for $q = c_{k}$, we have
  \begin{equation}
    \label{eq:RDHJ-RDHEq1}
    c_{k+1} = D_{2}H_{\rm d}^{+}\parentheses{ c_{k}, p_{k} },
  \end{equation}
  where we used Eq.~\eqref{eq:RDHJ-curve} and \eqref{eq:RDHJ_p-dS}.
  On the other hand, taking the derivative of Eq.~\eqref{eq:RDHJEq} with respect to $q$, 
  \begin{multline*}
    DS_{\rm d}^{k+1}(f^{+}_{k}(q)) \cdot Df^{+}_{k}(q) - DS_{\rm d}^{k}(q)
    - Df^{+}_{k}(q) \cdot D^{2}S_{\rm d}^{k+1}(f^{+}_{k}(q)) \cdot f^{+}_{k}(q)
    - DS_{\rm d}^{k+1}(f^{+}_{k}(q)) \cdot Df^{+}_{k}(q)
    \\
    + D_{1}H_{\rm d}^{+}\parentheses{ q, DS_{\rm d}^{k+1}(f^{+}_{k}(q)) }
    + D_{2}H_{\rm d}^{+}\parentheses{ q, DS_{\rm d}^{k+1}(f^{+}_{k}(q)) } \cdot D^{2}S_{\rm d}^{k+1}(f^{+}_{k}(q)) \cdot Df^{+}_{k}(q)
    = 0,
  \end{multline*}
  which reduces to 
  \begin{equation*}
    - DS_{\rm d}^{k}(q) + D_{1}H_{\rm d}^{+}\parentheses{ q, DS_{\rm d}^{k+1}(f^{+}_{k}(q)) } = 0,
  \end{equation*}
  due to Eq.~\eqref{eq:f_k-D2Hd+}.
  Then, substituting $q = c_{k}$ gives
  \begin{equation*}
    - DS_{\rm d}^{k}(c_{k}) + D_{1}H_{\rm d}^{+}\parentheses{ c_{k}, DS_{\rm d}^{k+1}(f^{+}_{k}(c_{k})) } = 0.
  \end{equation*}
  Using Eqs.~\eqref{eq:RDHJ-curve} and \eqref{eq:RDHJ_p-dS}, we obtain
  \begin{equation}
    \label{eq:RDHJ-RDHEq2}
    p_{k} = D_{1}H_{\rm d}^{+}\parentheses{ c_{k}, p_{k+1} }.
  \end{equation}
  Eqs.~\eqref{eq:RDHJ-RDHEq1} and \eqref{eq:RDHJ-RDHEq2} show that the sequence $\{(c_{k}, p_{k})\}$ satisfies the right discrete Hamilton's equations~\eqref{eq:RDHEq}.

  Now, let us prove the latter assertion.
  First, recall the implicit definition of $f^{-}_{k}$ in Eq.~\eqref{eq:f^-_k-implicit}:
  \begin{equation}
    \label{eq:f_k-D1Hd-}
    q = -D_{1}H_{\rm d}^{-}\parentheses{ DS_{\rm d}^{k}(q), f^{-}_{k}(q) }
  \end{equation}
  In particular, for $q = c_{k}$, we have
  \begin{equation}
    \label{eq:LDHJ-LDHEq1}
    c_{k} = -D_{1}H_{\rm d}^{-}\parentheses{ p_{k}, c_{k+1} },
  \end{equation}
  where we used Eq.~\eqref{eq:LDHJ-curve} and \eqref{eq:LDHJ_p-dS}.
  On the other hand, taking the derivative of Eq.~\eqref{eq:RDHJEq} with respect to $q$ yields, 
  \begin{multline*}
    DS_{\rm d}^{k+1}(f^{-}_{k}(q)) \cdot Df^{-}_{k}(q) - DS_{\rm d}^{k}(q)
    + D^{2}S_{\rm d}^{k}(q) \cdot q
    + DS_{\rm d}^{k}(q) 
    \\
    + D_{1}H_{\rm d}^{-}\parentheses{ DS_{\rm d}^{k}(q), f^{-}_{k}(q) } \cdot D^{2}S_{\rm d}^{k}(q)
    + D_{2}H_{\rm d}^{-}\parentheses{ DS_{\rm d}^{k}(q), f^{-}_{k}(q) } \cdot Df^{-}_{k}(q) = 0,
  \end{multline*}
  which reduces to
  \begin{equation*}
    \brackets{ 
      DS_{\rm d}^{k+1}(f^{-}_{k}(q))
    + D_{2}H_{\rm d}^{-}\parentheses{ DS_{\rm d}^{k}(q), f^{-}_{k}(q) }
    }
    \cdot Df^{-}_{k}(q) = 0,
  \end{equation*}
  due to Eq.~\eqref{eq:f_k-D1Hd-}.
  Then, substituting $q = c_{k}$ gives
  \begin{equation*}
    DS_{\rm d}^{k+1}(f^{-}_{k}(c_{k}))
    = -D_{2}H_{\rm d}^{-}\parentheses{ DS_{\rm d}^{k}(c_{k}), f^{-}_{k}(c_{k}) },
  \end{equation*}
  since $Df^{-}_{k}(c_{k})$ is invertible by assumption.
Then, using Eqs.~\eqref{eq:LDHJ-curve} and \eqref{eq:LDHJ_p-dS}, we obtain
  \begin{equation}
    \label{eq:LDHJ-LDHEq2}
    p_{k+1} = -D_{2}H_{\rm d}^{-}\parentheses{ p_{k}, c_{k+1} }.
  \end{equation}
  Eqs.~\eqref{eq:LDHJ-LDHEq1} and \eqref{eq:LDHJ-LDHEq2} show that the sequence $\{(c_{k}, p_{k})\}$ satisfies the left discrete Hamilton's equations~\eqref{eq:LDHEq}.
\end{proof}

\section{Application To Discrete Linear Hamiltonian Systems}
\label{sec:LinearHamiltonianSystems}
%In this section we take simple examples and solve the left discrete Hamilton--Jacobi equation~\eqref{eq:LDHJEq} for it.
%We expect to extend this to more complicated systems and also develop numerical methods based on the approach taken here.
\subsection{Discrete Linear Hamiltonian Systems and Matrix Riccati Equation}
\label{ssec:Riccati}
\begin{example}[Quadratic discrete Hamiltonian---discrete linear Hamiltonian systems]
  Consider a discrete Hamiltonian system on $T^{*}\R^{n} \cong \R^{n} \times \R^{n}$ (the configuration space is $Q = \R^{n}$) defined by the quadratic left discrete Hamiltonian
  \begin{equation}
    \label{eq:QDH}
    H_{\rm d}^{-}(p_{k}, q_{k+1}) = \frac{1}{2} p_{k}^{T} M^{-1} p_{k} + p_{k}^{T} L q_{k+1} + \frac{1}{2}q_{k+1}^{T} K q_{k+1},
  \end{equation}
  where $M$, $K$, and $L$ are real $n \times n$ matrices; we assume that $M$ and $L$ are invertible and also that $M$ and $K$ are symmetric.
  The left discrete Hamilton's equations~\eqref{eq:LDHEq} are
  \begin{equation*}
    \begin{array}{l}
      \DS q_{k} = -(M^{-1} p_{k} + L q_{k+1}),
      \medskip\\
      \DS p_{k+1} = -(L^{T} p_{k} + K q_{k+1}),
    \end{array}
  \end{equation*}
  or
  \begin{equation}
    \label{eq:DLHEq}
    \begin{pmatrix}
      q_{k+1} \medskip\\
      p_{k+1}
    \end{pmatrix}
    =
    \begin{pmatrix}
      -L^{-1} & -L^{-1} M^{-1} \medskip\\
      K L^{-1} & K L^{-1} M^{-1} - L^{T}
    \end{pmatrix}
    \begin{pmatrix}
      q_{k} \medskip\\
      p_{k}
    \end{pmatrix}.
  \end{equation}
  and hence are a discrete linear Hamiltonian system~(see Section~\ref{assec:LinearDiscreteHamiltonianSystems}).

  Now, let us solve the left discrete Hamilton--Jacobi equation~\eqref{eq:LDHJEq} for this system.
  For that purpose, we first generalize the problem to that with a set of initial points instead of a single initial point $(q_{0}, p_{0})$.
  More specifically, consider the set of initial points that is a Lagrangian affine space $\tilde{\mathcal{L}}_{z_{0}}$~(see Definition~\ref{def:LagrangianAffineSpace}) which contains the point $z_{0} \defeq (q_{0}, p_{0})$.
  Then, the dynamics is formally written as, for any discrete time $k \in \N$,
  \begin{equation*}
    \tilde{\mathcal{L}}^{(k)} \defeq (\tilde{F}_{L_{\rm d}})^{k}\parentheses{ \tilde{\mathcal{L}}_{z_{0}} } = \underbrace{\tilde{F}_{L_{\rm d}} \circ \dots \circ \tilde{F}_{L_{\rm d}}}_{k} \parentheses{ \tilde{\mathcal{L}}_{z_{0}} },
  \end{equation*}
  where $\tilde{F}_{L_{\rm d}}: T^{*}Q \to T^{*}Q$ is the discrete Hamiltonian map defined in Eq.~\eqref{eq:DHM}.
  Since $\tilde{F}_{L_{\rm d}}$ is a symplectic map, Proposition~\ref{prop:InvLagAffineSp} implies that $\tilde{\mathcal{L}}^{(k)}$ is a Lagrangian affine space.
  Then, assuming that $\tilde{\mathcal{L}}^{(k)}$ is transversal to $\{0\} \oplus Q^{*}$, Corollary~\ref{cor:GenFcn-LagAffineSp} implies that there exists a set of functions $S_{\rm d}^{k}$ of the form
  \begin{equation}
    \label{eq:S_d-QDH}
    S_{\rm d}^{k}(q) = \frac{1}{2} q^{T} A_{k} q + b_{k}^{T} q + c_{k},
  \end{equation}
  such that $\tilde{\mathcal{L}}^{(k)} = \Graph dS_{\rm d}^{k}$; here $A_{k}$ are symmetric $n \times n$ matrices, $b_{k}$ are elements in $\R^{n}$, and $c_{k}$ are in $\R$.

  Now that we know the form of the solution, we substitute the above expression into the discrete Hamilton--Jacobi equation to find the equations for $A_{k}$, $b_{k}$, and $c_{k}$.
  Notice first that the map $f^{-}_{k}$ is given by the first half of Eq.~\eqref{eq:DLHEq} with $p_{k}$ replaced by $DS_{\rm d}^{k}(q)$:
  \begin{align}
    f^{-}_{k}(q) &= -L^{-1}\parentheses{ q + M^{-1} DS_{\rm d}^{k}(q) }
    \nonumber\\
    &= -L^{-1}(I + M^{-1}A_{k})q - L^{-1}M^{-1}b_{k}.
    \label{eq:f^-_k-QDH}
  \end{align}
  Then, substituting Eq.~\eqref{eq:S_d-QDH} into the left-hand side of the left discrete Hamilton--Jacobi equation~\eqref{eq:LDHJEq} yields the following recurrence relations for $A_{k}$, $b_{k}$, and $c_{k}$:
  \begin{subequations}
    \label{eq:RR-QDH}
    \begin{align}
      A_{k+1} &= L^{T}(I + A_{k} M^{-1})^{-1} A_{k} L - K,
      \label{eq:RR-QDH-1}
      \\
      b_{k+1} &= -L^{T}(I + A_{k} M^{-1})^{-1} b_{k},
      \\
      c_{k+1} &= c_{k} - \frac{1}{2} b_{k}^{T} (M + A_{k})^{-1} b_{k},
    \end{align}
  \end{subequations}
  where we assumed that $I + A_{k} M^{-1}$ is invertible.
\end{example}

\begin{remark}
  For the $A_{k+1}$ defined by Eq.~\eqref{eq:RR-QDH-1} to be symmetric, it is sufficient that $A_{k}$ is invertible; for if it is, then Eq.~\eqref{eq:RR-QDH-1} becomes
  \begin{equation*}
    A_{k+1} = L^{T} (A_{k}^{-1} + M^{-1})^{-1} L - K,
  \end{equation*}
  and so $A_{k}$, $M$, and $K$ being symmetric implies that $A_{k+1}$ is as well.
\end{remark}

\begin{remark}
  We can rewrite Eq.~\eqref{eq:RR-QDH-1} as follows:
  \begin{equation}
    \label{eq:MREq}
    A_{k+1} = \brackets{ KL^{-1} + (K L^{-1} M^{-1} - L^{T}) A_{k} } (-L^{-1} - L^{-1}M^{-1}A_{k})^{-1}.
  \end{equation}
  Notice the exact correspondence between the coefficients in the above equation and the matrix entries in the discrete linear Hamiltonian equations~\eqref{eq:DLHEq}.
  In fact, this is the discrete Riccati equation that corresponds to the iteration defined by Eq.~\eqref{eq:DLHEq}.
  See \citet{AmMa1986} for details on this correspondence.
\end{remark}
To summarize the above observation, we have:
\begin{proposition}
  \label{prop:Riccati}
  The discrete Hamilton--Jacobi equation~\eqref{eq:LDHJEq} applied to the discrete linear Hamiltonian system~\eqref{eq:DLHEq} yields the discrete Riccati equation~\eqref{eq:MREq}.
\end{proposition}
In other words, the discrete Hamilton--Jacobi equation is a nonlinear generalization of the discrete Riccati equation.

\section{Relation to the Bellman Equation}
\label{sec:RelationToDTHJBEq}
In this section, we apply the above results to the optimal control setting.
We will show that the (right) discrete Hamilton--Jacobi equation \eqref{eq:RDHJEq} gives the Bellman equation (discrete-time HJB equation) as a special case.
This result gives a discrete analogue of the relationship between the H--J and HJB equations discussed in Section~\ref{sec:ConnectionWithOC}.

\subsection{Discrete Optimal Control Problem}
Let $q_{\rm d} \defeq \{ q_{k} \}_{k=0}^{N}$ be the state variables in a vector space $V \cong \R^{n}$ with $q_{0}$ and $q_{N}$ fixed and $u_{\rm d} \defeq \{ u_{k} \}_{k=0}^{N-1}$ be controls in the set $\mathcal{U} \subset \R^{m}$.
With a given function $C_{\rm d}: V \times \mathcal{U} \to \R$, define the discrete cost functional
\begin{equation*}
  J_{\rm d} \defeq \sum_{k=0}^{N-1} C_{\rm d}(q_{k}, u_{k}).
\end{equation*}
Then, we formulate the {\em Standard Discrete Optimal Control Problem} as follows~\citep[see, e.g.,][]{JoPo1964, Ca1970, Be1971, GuBl2004}:
\begin{problem}[Standard Discrete Optimal Control Problem]
  \label{prob:DOCP}
  Minimize the discrete cost functional, i.e., 
  \begin{equation}
    \label{eq:DOCP-optimality}
   \min_{u_{\rm d}} J_{\rm d}
   = 
   \min_{u_{\rm d}} \sum_{k=0}^{N-1} C_{\rm d}(q_{k}, u_{k}),
  \end{equation}
  subject to the constraint,
  \begin{equation}
    \label{eq:DOCP-constraint}
    q_{k+1} = f_{\rm d}(q_{k}, u_{k}).
  \end{equation}
\end{problem}

\subsection{Necessary Condition for Optimality and the Bellman Equation}
We would like to formulate the necessary condition for optimality.
First, introduce the augmented discrete cost functional:
\begin{align*}
  \hat{S}_{\rm d}(q_{\rm d}, p_{\rm d}, u_{\rm d})
  &\defeq \sum_{k=0}^{N-1} \braces{ C_{\rm d}(q_{k}, u_{k}) + p_{k+1} \cdot \brackets{ q_{k+1} - f_{\rm d}(q_{k}, u_{k})  } }
  \\
  &= \sum_{k=0}^{N-1} \brackets{ p_{k+1} \cdot q_{k+1} - \hat{H}_{\rm d}^{+}(q_{k}, p_{k+1}, u_{k}) },
\end{align*}
where we introduced the costate $p_{\rm d} \defeq \{ p_{k} \}_{k=1}^{N}$ with $p_{k} \in V^{*}$, and also defined the {\em discrete control Hamiltonian}
\begin{equation}
  \label{eq:Hhatd+}
  \hat{H}_{\rm d}^{+}(q_{k}, p_{k+1}, u_{k}) \defeq p_{k+1} \cdot f_{\rm d}(q_{k}, u_{k}) - C_{\rm d}(q_{k}, u_{k}).
\end{equation}
Then, the optimality condition,~Eq.~\eqref{eq:DOCP-optimality}, is restated as
\begin{equation*}
  \min_{q_{\rm d}, p_{\rm d}, u_{\rm d}} \hat{S}_{\rm d}(q_{\rm d}, p_{\rm d}, u_{\rm d})
  = \min_{q_{\rm d}, p_{\rm d}, u_{\rm d}} \sum_{k=0}^{N-1} \brackets{ p_{k+1} \cdot q_{k+1} - \hat{H}_{\rm d}^{+}(q_{k}, p_{k+1}, u_{k}) }.
\end{equation*}
In particular, extremality with respect to the control $u_{\rm d}$ implies
\begin{equation}
  D_{3}\hat{H}_{\rm d}^{+}(q_{k}, p_{k+1}, u_{k}) = 0,
  \quad
  k = 0, 1, \dots, N-1.
\end{equation}
Now, we assume that $\hat{H}_{\rm d}^{+}$ is sufficiently regular so that this equation uniquely determines the optimal control $u_{\rm d}^{*} \defeq \{ u_{k}^{*} \}_{k=0}^{N-1}$; and therefore, $u_{k}^{*}$ is a function of $q_{k}$ and $p_{k+1}$, i.e., $u_{k}^{*} = u_{k}^{*}(q_{k}, p_{k+1})$.
We then define
\begin{align}
  H_{\rm d}^{+}(q_{k}, p_{k+1})
  &\defeq \max_{u_{k}} \hat{H}_{\rm d}^{+}(q_{k}, p_{k+1}, u_{k})
  \nonumber\\
  &= \max_{u_{k}} \brackets{ p_{k+1} \cdot f_{\rm d}(q_{k}, u_{k}) - C_{\rm d}(q_{k}, u_{k}) }
  \nonumber\\
  &= p_{k+1} \cdot f_{\rm d}(q_{k}, u_{k}^{*}) - C_{\rm d}(q_{k}, u_{k}^{*}),
  \label{eq:RDH-OC}
\end{align}
and also the {\em optimal discrete cost-to-go function}
\begin{align}
  \label{eq:J_d^k}
  J_{\rm d}^{k}(q_{k}) &\defeq \sum_{l=k}^{N-1} \braces{ C_{\rm d}(q_{l}, u_{l}^{*}) + p_{l+1} \cdot \brackets{ q_{l+1} - f_{\rm d}(q_{l}, u_{l}^{*})  } }
  \nonumber\\
  &= \sum_{l=k}^{N-1} \brackets{ p_{l+1} \cdot q_{l+1} - H_{\rm d}^{+}(q_{l}, p_{l+1}) }
  \nonumber\\
  &= S_{\rm d}^{*} - S_{\rm d}^{k}(q_{k}),
\end{align}
where $S_{\rm d}^{*}$ is the {\em optimal discrete cost functional}, i.e., 
\begin{equation*}
  S_{\rm d}^{*} \defeq \hat{S}_{\rm d}(q_{\rm d}, p_{\rm d}, u_{\rm d}^{*})
  = \sum_{k=0}^{N-1} \brackets{ p_{k+1} \cdot q_{k+1} - H_{\rm d}^{+}(q_{k}, p_{k+1}) }.
\end{equation*}
and
\begin{equation*}
  S_{\rm d}^{k}(q_{k})
  \defeq \sum_{l=0}^{k-1} \brackets{ p_{l+1} \cdot q_{l+1} - H_{\rm d}^{+}(q_{l}, p_{l+1}) }.
\end{equation*}
The above action sum has exactly the same form as Eq.~\eqref{eq:S_d+} formulated in the framework of discrete Hamiltonian mechanics.
Therefore, our theory now directly applies to this case:
The corresponding right discrete Hamilton's equations~\eqref{eq:RDHEq} are, using the expression in Eq.~\eqref{eq:RDH-OC}, 
\begin{equation*}
  \begin{array}{l}
    \DS q_{k+1} = f_{\rm d}(q_{k}, u_{k}^{*}),
    \medskip\\
    \DS p_{k} = p_{k+1} \cdot D_{1}f_{\rm d}(q_{k}, u_{k}^{*}) - D_{1}C_{\rm d}(q_{k}, u_{k}^{*}).
  \end{array}
\end{equation*}
Therefore, Eq.~\eqref{eq:f^+_k-implicit} gives the implicit definition of $f^{+}_{k}$ as follows:
\begin{equation}
  \label{eq:f_k-OC}
  f^{+}_{k}(q_{k}) = f_{\rm d}\parentheses{ q_{k}, u_{k}^{*}\parentheses{ q_{k}, DS_{\rm d}^{k+1}(f^{+}_{k}(q_{k})) } }.
\end{equation}
Hence, the (right) discrete Hamilton--Jacobi equation~\eqref{eq:RDHJEq} applied to this case gives
\begin{equation*}
  S_{\rm d}^{k+1}(f_{\rm d}(q_{k}, u_{k}^{*})) - S_{\rm d}^{k}(q_{k}) 
  - DS_{\rm d}^{k+1}(f_{\rm d}(q_{k}, u_{k}^{*})) \cdot f_{\rm d}(q_{k}, u_{k}^{*})
  + H_{\rm d}^{+}\parentheses{ q_{k}, DS_{\rm d}^{k+1}(f_{\rm d}(q_{k}, u_{k}^{*})) }
  = 0,
\end{equation*}
and again using the expression for the Hamiltonian in Eq.~\eqref{eq:RDH-OC}, this becomes
\begin{equation*}
  \max_{u_{k}} \brackets{
    S_{\rm d}^{k+1}(f_{\rm d}(q_{k}, u_{k})) - C_{\rm d}(q_{k}, u_{k})
  }
  = S_{\rm d}^{k}(q_{k}).
\end{equation*}
Since $S_{\rm d}^{k}(q_{k}) = S_{\rm d}^{*} - J_{\rm d}^{k}(q_{k})$, we obtain
\begin{equation}
  \label{eq:DHJBEq}
  \min_{u_{k}}\brackets{
    J_{\rm d}^{k+1}(f_{\rm d}(q_{k}, u_{k})) + C_{\rm d}(q_{k}, u_{k})
  }
  = J_{\rm d}^{k}(q_{k}),
\end{equation}
which is the {\em Bellman equation}~(see, e.g., \citet{Be1971,Be1972} and \citet{Be2005}).

\begin{remark}
  \label{remark:simplicity_of_DHJBEq}
  Notice that the discrete HJB equation~\eqref{eq:DHJBEq} is much simpler than the discrete Hamilton--Jacobi equations~\eqref{eq:RDHJEq} and \eqref{eq:LDHJEq} because of the special form of the control Hamiltonian~Eq.~\eqref{eq:RDH-OC}.
  Also, notice that, as shown in Eq.~\eqref{eq:f_k-OC}, the term $f^{+}_{k}(q_{k})$ is written in terms of the given function $f$.
  See Remark~\ref{remark:complexity_of_DHJEq} for comparison.
\end{remark}

\subsection{Relation between the Discrete H--J and Bellman Equations and its Consequences}
Summarizing the observation made above, we have
\begin{proposition}
  \label{prop:RDHJ-DHJB}
  The right discrete Hamilton--Jacobi equation~\eqref{eq:RDHJEq} applied to the Hamiltonian formulation of the Standard Discrete Optimal Control Problem~\ref{prob:DOCP} gives the Bellman equation~\eqref{eq:DHJBEq}.
\end{proposition}

This observation leads to the following well-known fact:
\begin{proposition}
  Let $J_{\rm d}^{k}(q_{k})$ be a solution to the Bellman equation~\eqref{eq:DHJBEq}.
  Then, the costate $p_{k}$ in the discrete maximum principle is given as follows:
  \begin{equation*}
    p_{k} = -DJ_{\rm d}^{k}(c_{k}),
  \end{equation*}
  where $c_{k+1} = f_{\rm d}(c_{k}, u_{k}^{*})$ with the optimal control $u_{k}^{*}$.
\end{proposition}

\begin{proof}
  Follows from a reinterpretation of Theorem~\ref{thm:DHJ} through Proposition~\ref{prop:RDHJ-DHJB} with the relation $S_{\rm d}^{k}(q_{k}) = S_{\rm d}^{*} - J_{\rm d}^{k}(q_{k})$.
\end{proof}

\section{Generalized Bellman Equation with Internal-Stage Controls}
\label{sec:Higher-Order_Bellman}
In the previous section, we showed that the discrete Hamilton--Jacobi equation recovers the Bellman equation if we apply our theory to the Hamiltonian formulation of the Standard Discrete Optimal Control Problem~\ref{prob:DOCP}.
In this section, we generalize the approach to derive what may be considered as higher-order discrete-time approximations of the HJB equation~\eqref{eq:HJBEq}.
Namely, we derive a class of discrete control Hamiltonians that use higher-order approximations (a more general version of Eq.~\eqref{eq:Hhatd+}) by employing the technique of Galerkin Hamiltonian variational integrators introduced by \citet{LeZh2011}; and then, we apply our theory to obtain a class of generalized Bellman equations that have controls at internal stages.

\subsection{Continuous-Time Optimal Control Problem}
Let us first briefly review the standard formulation of continuous-time optimal control problems.
Let $q$ be the state variable in a vector space $V \cong \R^{n}$, $q_{0}$ and $q_{T}$ fixed in $V$, and $u$ be the control in the set $\mathcal{U} \subset \R^{m}$.
With a given function $C: V \times \mathcal{U} \to \R$, define the cost functional
\begin{equation*}
  J \defeq \int_{0}^{T} C(q(t), u(t))\,dt.
\end{equation*}
Then, we formulate the {\em Standard Continuous-Time Optimal Control Problem} as follows:
\begin{problem}[Standard Continuous-Time Optimal Control Problem]
  \label{prob:OCP}
  Minimize the cost functional, i.e., 
  \begin{equation*}
    \min_{u(\cdot)} J
    = 
    \min_{u(\cdot)} \int_{0}^{T} C(q(t), u(t))\,dt,
  \end{equation*}
  subject to the constraints,
  \begin{equation*}
    \dot{q} = f(q, u),
  \end{equation*}
  and $q(0) = q_{0}$ and $q(T) = q_{T}$.
\end{problem}

A Hamiltonian structure comes into play with the introduction of the augmented cost functional:
\begin{align*}
  \hat{S} &\defeq \int_{0}^{T} \braces{ C(q(t), u(t)) + p(t) \brackets{ \dot{q}(t) - f(q(t), u(t)) } } dt
  \\
  &= \int_{0}^{T} \brackets{ p(t) \dot{q}(t) - \hat{H}(q(t), p(t), u(t)) } dt,
\end{align*}
where we introduced the costate $p(t) \in V^{*}$, and also defined the {\em control Hamiltonian},
\begin{equation}
  \label{eq:Hhat}
  \hat{H}(q, p, u) \defeq p \cdot f(q, u) - C(q, u).
\end{equation}

\subsection{Galerkin Hamiltonian Variational Integrator}
Recall, from \citet[Section~2.2]{LeZh2011}, that the exact right discrete Hamiltonian is a type-two generating function for the original continuous-time Hamiltonian flow, defined by
\begin{equation}
  \label{eq:Hdex+}
  H_{\rm d, ex}^{+}(q_0, p_{1})
  = \ext_{\substack{(q, p) \in \mathcal{C}^1([0, h],T^*Q)\\q(0)=q_0, p(h)=p_1}} \braces{
    p_1 q_1 - \int_0^h \left[ p(t) \dot{q}(t)-H(q(t), p(t)) \right]dt
  },
\end{equation}
where $h$ is the time step; $\mathcal{C}^1([0, h],T^*Q)$ is the set of continuously differentiable curves on $T^{*}Q$ over the time interval $[0, h]$; an extremum is achieved for the exact solution of Hamilton's equations~\eqref{eq:HamiltonsEq} that satisfy the specified boundary conditions.
Therefore, it requires the exact solution $(q(t), p(t))$ to evaluate the the above integral, and so the exact discrete Hamiltonian cannot be practically computed in general.

The key idea of Galerkin Hamiltonian variational integrators~\cite{LeZh2011} is to replace the set of curves $\mathcal{C}^1([0, h],T^*Q)$ by a certain finite-dimensional space so as to obtain a computable expression for a discrete Hamiltonian.

\subsection{Galerkin Discrete Control Hamiltonian}
Here, we would like to apply the above idea to the {\em control} Hamiltonian, Eq.~\eqref{eq:Hhat}, to obtain a discrete control Hamiltonian.

Let $\mathcal{C}_{\rm d}^s(V)$ be a finite-dimensional space of curves defined by
\begin{equation*}
  \mathcal{C}_{\rm d}^s(V) \defeq \setdef{ t \mapsto \sum_{i=1}^{s} w^{i}\,\psi_{i}(t/h) }{ t \in [0,h],\ w^{i} \in V \text{ for each } i \in \{1, 2, \dots, s\} }.
\end{equation*}
with the basis functions $\{ \psi_{i}: [0,1] \to \R \}_{i=1}^{s}$.
\begin{enumerate}[\sf 1.]
\item Use the basis functions $\psi_i$ to approximate the velocity $\dot{q}$ over the interval $[0,h]$,
  \[
  \dot{q}(\tau h) \approx \dot{q}_{\rm d}(\tau h) = \sum_{i=1}^s w^i \psi_i(\tau),
  \]
  where $\tau \in [0,1]$ and $w^{i} \in V$ for each $i = 1, \dots, s$.
  \medskip
\item Integrate $\dot{q}_{\rm d}(t)$ over $[0, \tau h],$ to obtain the approximation for the position $q$, i.e.,
  \[
  q_{\rm d}(\tau h)
  = q_{\rm d}(0)+ \int_0^{\tau h } \sum_{i=1}^s w^i \psi_i(t/h)\,dt
  = q_0 + h\sum_{i=1}^s w^i \int_0^{\tau} \psi_i(\rho)\,d\rho,
  \]
  where we applied the boundary condition $q_{\rm d}(0)=q_0$.
  Applying the boundary condition $q_{\rm d}(h)=q_1$ at the other endpoint yields
  \[
  q_1
  = q_{\rm d}(h)
  = q_{0} + h\sum_{i=1}^{s} w^{i} \int_{0}^{1} \psi_{i}(\rho)\,d\rho
  = q_{0} + h\sum_{i=1}^{s} B_{i} w^{i},
  \]
  where $B_i \defeq \int_0^{1} \psi_i(\tau)\,d\tau$.
  Furthermore, we introduce the internal stages,
  \begin{equation}
    \label{eq:Q}
    Q^{i}(w) \defeq q_{\rm d}(c^{i} h)
    = q_0 + h \sum_{j=1}^{s} w^{j} \int_0^{c^{i}} \psi_{j}(\tau)\,d\tau
    = q_0 + h \sum_{j=1}^{s} A^{i}_{j }w^{j},
  \end{equation}
  where $A^{i}_{j} \defeq \int_0^{c^i} \psi_j(\tau)\, d\tau$.
  \medskip
\item The exact discrete control Hamiltonian $\hat{H}_{\rm d, ex}^{+}$ is defined as in Eq.~\eqref{eq:Hdex+}:
  \begin{equation*}
    \hat{H}_{\rm d, ex}^{+}(q_0, p_{1}; u(\cdot)) \defeq
    \ext_{\substack{(q, p) \in \mathcal{C}^1([0, h],V \times V^*)\\q(t_0)=q_0,\, p(t_1)=p_{1}}} \braces{
    p_{1} q_1 - \int_0^h \left[ p(t) \dot{q}(t) - \hat{H}(q(t), p(t), u(t)) \right] dt
    }.
  \end{equation*}
  Again this is practically not computable, and so we employ the following approximation:
  Use the numerical quadrature formula
  \begin{equation*}
    \int_{0}^{1}f(\rho)\,d\rho \approx \sum_{i=1}^{s} b_{i}\, f(c^{i})
  \end{equation*}
  with constants $(b_i, c^i)$ and the finite-dimensional function space $\mathcal{C}_{\rm d}^s(V)$ to construct $\hat{H}_{\rm d}^+(q_0, p_{1}, U)$ as follows:
  \begin{align*}
    \hat{H}_{\rm d}^+(q_0, p_{1}, U)
    &\defeq \ext_{\substack{\dot{q}_{\rm d}\in \mathcal{C}_{\rm d}^s(V)\\P^{i}\in V^*}} \left\{
      p_{1} q_{\rm d}(h)
      - h\sum_{i=1}^s b_i\left[
        p(c^i h) \dot{q}_{\rm d}(c^i h)
        - \hat{H}\parentheses{ Q^{i}(w), P^{i}, U^{i} }
      \right]
    \right\}
    \nonumber\\
    &= \ext_{w,P} K(q_0, w, P, U, p_{1}),
  \end{align*}
  where we set $P^i \defeq p(c^i h)$ and $U^{i} \defeq u(c^{i} h)$ and defined
  \begin{align*}
    K(q_0, w, P, U, p_{1})
    &\defeq p_{1} \cdot \left( q_0 + h \sum_{i=1}^{s} B_{i} w^{i} \right)
    -h\sum_{i=1}^s b_i \left[
      P^i \cdot \sum_{j=1}^s w^j \psi_j(c^i) - \hat{H}\left( Q^{i}(w), P^{i}, U^{i}\right)
    \right]
    \\
    &= p_{1} \cdot \left( q_0 + h \sum_{i=1}^{s} B_{i} w^{i} \right)
    \\
    &\quad- h\sum_{i=1}^s b_i \left\{
      P^i \cdot \brackets{ \sum_{j=1}^{s} M^{i}_{j} w^{j} - f(Q^{i}(w), U^{i}) } + C(Q^{i}(w), U^{i})
    \right\},
   \end{align*}
  where we defined $M^{i}_{j} \defeq \psi_{j}(c^{i})$ and used the expression for the control Hamiltonian in Eq.~\eqref{eq:Hhat}; note that $P^{i} \in V^{*}$ and $w^{i} \in V$ for each $i = 1, \dots, s$, and that $f$ takes values in $V$.
  In order to obtain an expression for $H_{\rm d}^+(q_0, p_{1}, U)$, we first compute the stationarity conditions for $K(q_0, w, P, U, p_{1})$ under the fixed boundary condition $(q_0,p_{1})$:
  \begin{subequations}
    \label{gsta}
    \begin{align}
      0 &= \frac{\partial K(q_0, w, P, U, p_{1})}{\partial w^j}
      = h p_{1} \cdot B_j
      - h\sum_{i=1}^{s} b_i \left[
        M^{i}_{j}\, P^{i}
        - h A^{i}_{j}\, D_{1}\hat{H}\parentheses{ q^i(w),P^i,U^{i} }
      \right],
      \label{gsta1}
      \\
      0 &= \frac{\partial  K(q_0, w, P, U, p_{1})}{\partial P^i}
      = -h b_i
      \left[
        \sum_{j=1}^{s} M^{i}_{j}\, w^{j} - f(Q^{i}(w), U^{i})
      \right],
      \label{gsta2}
    \end{align}
  \end{subequations}
  for $j = 1, \dots, s$.
  \medskip
\item By solving the $2 s$ stationarity conditions \eqref{gsta}, we can express the parameters $w$ and $P$ in terms of $q_0$, $p_{1}$, and $U$, i.e., $w = \tilde{w}(q_0, U)$ and $P = \tilde{P}(q_0, p_{1}, U)$: In particular, assuming $b_{j} \neq 0$ for each $j = 1, \dots, s$, Eq.~\eqref{gsta2} gives $w^{j} M_{j}^{i} = f(Q^{i}(w), U^{i})$; this gives a set of $n s$ nonlinear equations\footnote{Recall that $w^{j} \in V$ and $f$ takes values in $V$.} satisfied by $w = \tilde{w}(q_0, U)$.\footnote{Note from Eq.~\eqref{eq:Q} that $Q^{i}$ is written in terms of $q_{0}$ and $w$.}
  Therefore, we have
  \begin{equation*}
    K\parentheses{ q_0, \tilde{w}(q_0, U), P, U, p_{1} }
    = p_{1} \cdot \left[ q_0 + h \sum_{j=1}^{s} B_{j}\, \tilde{w}^{j}(q_0, U) \right]
    - h\sum_{i=1}^s b_i\,C\parentheses{ Q^{i}(\tilde{w}(q_0, U)), U^{i} }.
  \end{equation*}
  Notice that the internal-stage momenta, $P^{i}$, disappear when we substitute $w = \tilde{w}(q_{0}, U)$.
  Therefore, we obtain the following Galerkin discrete control Hamiltonian:
  \begin{align}
    \label{eq:Hhatd+-Galerkin}
    \hat{H}_{\rm d}^{+}(q_{0}, p_{1}, U) &\defeq
    K\parentheses{ q_0, \tilde{w}(q_0, U), \tilde{P}(q_0, p_{1}, U), U, p_{1} }
    \nonumber\\
    &= p_{1} \cdot \left[ q_0 + h \sum_{j=1}^{s} B_{j}\, \tilde{w}^{j}(q_0, U) \right]
    - h \sum_{i=1}^s b_i\,C\parentheses{ Q^{i}(\tilde{w}(q_0, U)), U^{i} }.
  \end{align}
%  Then we obtain
%  \begin{align}
%    \label{gdefq1}
%    q_1 &= \frac{\partial}{\partial p_{1}} H_{\rm d}^+(q_0, p_{1}, U)
%    = \frac{\partial}{\partial p_{1}} K\parentheses{ q_0, \tilde{w}(q_0, U), \tilde{P}(q_0, p_{1}, U), U, p_{1} }
%    \nonumber\\
%    &= \sum_{i=1}^{s} \frac{\partial K}{\partial P^i} \cdot \frac{\partial \tilde{P}^i}{\partial p_{1}}
%    + \frac{\partial K}{\partial p_{1}}
%    = \frac{\partial K}{\partial p_{1}}
%    \nonumber\\
%    &= q_0 + h\sum_{i=1}^s  B_i\, \tilde{w}^i(q_{0}, U).
%  \end{align}
\end{enumerate}

\subsection{The Bellman Equation with Internal-Stage Controls}
\begin{figure}[ht!]
  \centering
  \includegraphics[width=.5\linewidth]{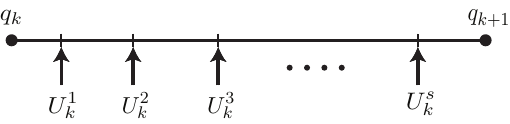}
  \caption{Internal-stage controls $\{U_{k}^{i}\}_{i=1}^{s}$ in the discrete time intervals between $k$ and $k+1$.}
  \label{fig:InternalStageControls}
\end{figure}

The Galerkin discrete control Hamiltonian, Eq.~\eqref{eq:Hhatd+-Galerkin}, gives
\begin{equation}
  \label{eq:Hhatd+-Galerkin2}
  \hat{H}_{\rm d}^{+}(q_{k}, p_{k+1}, U_{k}^{1}, \dots, U_{k}^{s}) \defeq
  p_{k+1} \cdot f_{\rm d}(q_{k}, U_{k}^{1}, \dots, U_{k}^{s}) - C_{\rm d}(q_{k}, U_{k}^{1}, \dots, U_{k}^{s}),
\end{equation}
with
\begin{equation}
  \label{eq:f_d-Galerkin}
  f_{\rm d}(q_{k}, U_{k}^{1}, \dots, U_{k}^{s}) \defeq
  q_k + h\sum_{i=1}^s  \tilde{w}^i(q_{k}, U_{k}^{1}, \dots, U_{k}^{s})\,B_i,
\end{equation}
and
\begin{equation}
  \label{eq:Cd-Galerkin}
  C_{\rm d}(q_{k}, U_{k}^{1}, \dots, U_{k}^{s}) \defeq h \sum_{i=1}^s b_i\,C\parentheses{ Q^{i}(\tilde{w}(q_k, U_{k}^{1}, \dots, U_{k}^{s})), U_{k}^{i} }.
\end{equation}
This is a generalized version of Eq.~\eqref{eq:Hhatd+} with internal-stage controls $\{ U_{k}^{i} \}_{i=1}^{s}$ as opposed to a single control $u_{k}$ per time step (see Fig.~\ref{fig:InternalStageControls}).
Now assume that
\begin{equation}
  \pd{}{U_{k}^{i}}\hat{H}_{\rm d}^{+}(q_{k}, p_{k+1}, U_{k}^{1}, \dots, U_{k}^{s}) = 0,
  \quad
  i = 1, \dots, s,
\end{equation}
is solvable for $\{ U_{k}^{i} \}_{i=1}^{s}$ to give the optimal internal-stage controls $\{ U_{k}^{*,i} \}_{i=1}^{s}$.
Then, we may apply the same argument as in Section~\ref{sec:RelationToDTHJBEq}: In particular, the right discrete Hamilton--Jacobi equation~\eqref{eq:RDHJEq} applied to this case gives the following {\em Bellman equation with internal-stage controls}:
\begin{equation}
  \label{eq:DHJBEq-Galerkin}
  \min_{U_{k}^{1},\dots,U_{k}^{s}}\brackets{ J_{\rm d}^{k+1}(f_{\rm d}(q_{k}, U_{k}^{1}, \dots, U_{k}^{s})) + C_{\rm d}(q_{k}, U_{k}^{1}, \dots, U_{k}^{s}) }
  = J_{\rm d}^{k}(q_{k}).
\end{equation}
%the minimization problem
%\begin{equation}
%  \label{eq:DOCP-optimality-Galerkin}
%  \min_{u_{\rm d}} \sum_{k=0}^{N-1} C_{\rm d}(q_{k}, U_{k}^{1}, \dots, U_{k}^{s})
%\end{equation}
%subject to the constraint
%\begin{equation}
%  \label{eq:DOCP-constraint-Galerkin}
%  q_{k+1} = f_{\rm d}(q_{k}, U_{k}^{1}, \dots, U_{k}^{s}),
%\end{equation}

The following example shows that the standard Bellman equation~\eqref{eq:DHJBEq} follows as a special case:
\begin{example}[The standard Bellman equation]
  \label{ex:StandardBellmanEq}
  Let $s = 1$, and select
  \begin{equation*}
    \Psi_{1}(\tau) = 1,
    \qquad
    b_{1} = 1,
    \qquad
    c^{1} = 0.
  \end{equation*}
  Then, we have $B_{1} = 1$, $A^{1}_{1} = 0$, and $M^{1}_{1} = 1$.
  Hence, Eq.~\eqref{eq:Q} gives $Q^{1} = q_{k}$ (we set the endpoints $(q_{0}, q_{1})$ to be $(q_{k}, q_{k+1})$ here), and Eq.~\eqref{gsta2} gives
  \begin{equation*}
    w^{1} = f(q_{k}, U_{k}^{1}).
  \end{equation*}
 However, the control $U_{k}^{1}$ is defined as follows (we shift the time intervals from $[0,h]$ to $[t_{k}, t_{k} + h]$ here):
 \begin{equation*}
   U_{k}^{1} \defeq u(t_{k} + c^{1}h) = u(t_{k}) \eqdef u_{k}.
 \end{equation*}
 So, we have $w^{1} = f(q_{k}, u_{k})$, and thus, Eqs.~\eqref{eq:f_d-Galerkin} and \eqref{eq:Cd-Galerkin} give
 \begin{equation*}
   f_{\rm d}(q_{k}, u_{k}) = q_{k} + h\,f(q_{k}, u_{k}),
 \end{equation*}
 and
 \begin{equation*}
   C_{\rm d}(q_{k}, u_{k}) = h\,C(q_{k}, u_{k}),
 \end{equation*}
 respectively.
 Notice that this approximation gives the forward-Euler discretization of the Standard Continuous-Time Optimal Control Problem~\ref{prob:OCP} to yield the Standard Discrete Optimal Control Problem~\ref{prob:DOCP}.
 In fact, the Bellman equation with internal-stage controls, Eq.~\eqref{eq:DHJBEq-Galerkin}, reduces to the standard Bellman equation~\eqref{eq:DHJBEq}:
 \begin{equation}
   \label{eq:StandardBellmanEq}
   \min_{u_{k}}\brackets{
     J_{\rm d}^{k+1}(f_{\rm d}(q_{k}, u_{k})) + C_{\rm d}(q_{k}, u_{k})
   }
   = J_{\rm d}^{k}(q_{k}).
 \end{equation}
\end{example}

\begin{remark}
  Higher-order approximations with any number of $s$ are possible as long as Eq.~\eqref{gsta2} is solvable for $w$.
  See \citet{LeZh2011} for various different choices of discretizations.
\end{remark}

\subsection{Application to the Heisenberg System}
Let us now apply the above results to a simple optimal control problem to illustrate the result:
\begin{example}[The Heisenberg system; see, e.g., \citet{Br1981} and \citet{Bl2003}]
  Consider the following optimal control problem: For a fixed time $T > 0$,
  \begin{equation*}
    \min_{u(\cdot),\,v(\cdot)} \int_{0}^{T} \frac{1}{2}(u^{2} + v^{2})\,dt,
  \end{equation*}
  subject to the constraint,
  \begin{equation*}
    \dot{x} = u, 
    \qquad
    \dot{y} = v,
    \qquad
    \dot{z} = u y - v x.
  \end{equation*}
  This is the Standard Continuous-Time Optimal Control Problem~\ref{prob:OCP} with $V = \R^{3}$, $\mathcal{U} = \R^{2}$, $q = (x, y, z)$, and
  \begin{equation*}
    f(x, y, z, u, v) =
    \begin{pmatrix}
      u \\
      v \\
      u y - v x
    \end{pmatrix},
    \qquad
    C(x, y, z, u, v) = \frac{1}{2}(u^{2} + v^{2}).
  \end{equation*}
  If we apply the choice of the discretization in Example~\ref{ex:StandardBellmanEq}, we have the standard Bellman equation
  \begin{equation*}
    \min_{u_{k}, v_{k}}\brackets{ S_{\rm d}^{k+1}(f_{\rm d}(q_{k}, u_{k}, v_{k})) - C_{\rm d}(q_{k}, u_{k}, v_{k}) }
    - S_{\rm d}^{k}(q_{k}) 
    = 0,
  \end{equation*}
  with $q_{k} \defeq (x_{k}, y_{k}, z_{k})$, where
  \begin{equation*}
    f_{\rm d}(q_{k}, u_{k}, v_{k}) =
    \begin{pmatrix}
      x_{k} + h\,u_{k} \\
      y_{k} + h\,v_{k} \\
      z_{k} + h\,(u_{k} y_{k} - v_{k} x_{k})
    \end{pmatrix},
    \qquad
    C_{\rm d}(q_{k}, u_{k}, v_{k}) = \frac{h}{2}(u_{k}^{2} + v_{k}^{2}).
  \end{equation*}

  Now, if we choose $s = 2$, and select
  \begin{equation*}
    (\Psi_{1}(\tau), \Psi_{2}(\tau)) = (1, \cos(\pi\tau)),
    \qquad
    b = (b_{1}, b_{2}) = \parentheses{ \frac{1}{2}, \frac{1}{2} },
    \qquad
    c = (c^{1}, c^{2}) = (0, 1).
  \end{equation*}
  Then, we have
  \begin{equation*}
    B = (B_{1}, B_{2}) = (1, 0),
    \qquad
    A =
    \begin{pmatrix}
      0 & 0 \\
      1 & 0
    \end{pmatrix},
    \qquad
    M = 
    \begin{pmatrix}
      1 & 1 \\
      1 & -1
    \end{pmatrix}.
  \end{equation*}
  \citet[Example~4.4]{LeZh2011} show that this choice of discretization corresponds to the St\"ormer--Verlet method~(see, e.g., \citet{MaWe2001}).
  The Bellman equation with internal-stage controls, Eq.~\eqref{eq:DHJBEq-Galerkin}, then becomes
  \begin{equation*}
    \min_{ u_{k}^{1},v_{k}^{1},u_{k}^{2},v_{k}^{2} }\brackets{
      J_{\rm d}^{k+1}\parentheses{ f_{\rm d}\parentheses{ q_{k}, u_{k}^{1}, v_{k}^{1}, u_{k}^{2}, v_{k}^{2} } }
      + C_{\rm d}\parentheses{ q_{k}, u_{k}^{1}, v_{k}^{1}, u_{k}^{2}, v_{k}^{2} }
    }
    = J_{\rm d}^{k}(q_{k}),
  \end{equation*}
  where
  \begin{equation*}
    f_{\rm d}\parentheses{ q_{k}, u_{k}^{1}, v_{k}^{1}, u_{k}^{2}, v_{k}^{2} } =
    \begin{pmatrix}
      x_{k} + h\,(u_{k}^{1} + u_{k}^{2})/2
      \medskip\\
      y_{k} + h\,(v_{k}^{1} + v_{k}^{2})/2
      \medskip\\
      \DS z_{k} + h \parentheses{
        \frac{u_{k}^{1} + u_{k}^{2}}{2}\,y_{k} - \frac{v_{k}^{1} + v_{k}^{2}}{2}\, x_{k}
        + \frac{u_{k}^{2} v_{k}^{1} - u_{k}^{1} v_{k}^{2}}{4}
      }
    \end{pmatrix},
  \end{equation*}
  and
  \begin{equation*}
    C_{\rm d}\parentheses{ q_{k}, u_{k}^{1}, v_{k}^{1}, u_{k}^{2}, v_{k}^{2} }
    = \frac{h}{2}\brackets{ \frac{(u_{k}^{1})^{2} + (u_{k}^{2})^{2}}{2} + \frac{(v_{k}^{1})^{2} + (v_{k}^{2})^{2}}{2} }.
  \end{equation*}
\end{example}
%Let $s = 1$ and select
%\begin{equation*}
%  \text{(i)}~\parentheses{ c^{1}, b_{1}, A^{1}_{1} } = \parentheses{ 0, 1, 0 };
%  \qquad
%  \text{(ii)}~\parentheses{ c^{1}, b_{1}, A^{1}_{1} } = \parentheses{ \frac{1}{2}, 1, \frac{1}{2} };
%  \qquad
%  \text{(iii)}~\parentheses{ c^{1}, b_{1}, A^{1}_{1} } = \parentheses{ 1, 1, 1 }.
%\end{equation*}
%Then Eq.~\eqref{eq:Q} gives
%\begin{equation*}
%  \text{(i)}~Q^{1} = q_{k};
%  \qquad
%  \text{(ii)}~Q^{1} = q_{k} + \frac{1}{2}\,h\,w^{1};
%  \qquad
%  \text{(iii)}~Q^{1} = q_{k} + h\,w^{1}
%\end{equation*}

\section{Conclusion and Future Work}
We developed a discrete-time analogue of Hamilton--Jacobi theory starting from the discrete variational Hamiltonian mechanics formulated by \citet{LaWe2006}.
We reinterpreted and extended the discrete Hamilton--Jacobi equation given by \citet{ElSc1996} in the language of discrete mechanics.
Furthermore, we showed that the discrete Hamilton--Jacobi equation reduces to the discrete Riccati equation with a quadratic Hamiltonian, and also that it specializes to the Bellman equation of dynamic programming if applied to standard discrete optimal control problems.
These results are discrete analogues of the corresponding known results in the continuous-time theory.
Application to discrete optimal control also revealed that the Discrete Hamilton--Jacobi Theorem~\ref{thm:DHJ} specializes to a well-known result in discrete optimal control theory.
We also used a Galerkin-type approximation to derive Galerkin discrete control Hamiltonians.
This technique gave an explicit formula for discrete control Hamiltonians in terms of the constructs in the original continuous-time optimal control problem.
By viewing the Bellman equation as a special case of the discrete Hamilton--Jacobi equation, we could introduce the discretization technique for discrete Hamiltonian mechanics into the discrete optimal control setting; this lead us to a class of Bellman equations with controls at internal stages.

We are interested in the following topics for future work:
\begin{itemize}
\item {\em Application to integrable discrete systems}:
  Theorem~\ref{thm:DHJ} gives a discrete analogue of the theory behind the technique of solution by separation of variables, i.e., the theorem relates a solution of the discrete Hamilton--Jacobi equations with that of the discrete Hamilton's equations.
  An interesting question then is whether or not separation of variables applies to integrable discrete systems, e.g., discrete rigid bodies of \citet{MoVe1991} and various others discussed by \citet{Su2003,Su2004}.
  \smallskip
\item {\em Development of numerical methods based on the discrete Hamilton--Jacobi equation}:
  Hamilton--Jacobi equation has been used to develop structured integrators for Hamiltonian systems.
  \citet{GeMa1988} developed a numerical method that preserves momentum maps and Poisson brackets of Lie--Poisson systems by solving the Lie--Poisson Hamilton--Jacobi equation approximately.
  See also \citet{ChSc1990} (and references therein) for a survey of structured integrators based on the Hamilton--Jacobi equation.
  The present theory, being inherently discrete in time, potentially provides a variant of such numerical methods.
  \smallskip
\item {\em Extension to discrete nonholonomic and Dirac mechanics}:
  The present work is concerned only with unconstrained systems.
  Extensions to nonholonomic and Dirac mechanics, more specifically discrete-time versions of the nonholonomic Hamilton--Jacobi theory~\cite{IgLeMa2008, LeMaMa2010, OhBl2009, CaGrMaMaMuRo2010} and Dirac Hamilton--Jacobi theory~\cite{DiracHJ}, are another direction for future research.
  \smallskip
\item {\em Relation to the power method and iterations on the Grassmannian manifold}:
  \citet{AmMa1986} established links between the power method, iterations on the Grassmannian manifold, and the Riccati equation.
  The discussion on iterations of Lagrangian subspaces and its relation to the Riccati equation in Sections~\ref{ssec:Riccati} and \ref{ssec:LagrangianSubAffinceSpaces} is a special case of such links.
  On the other hand, Proposition~\ref{prop:Riccati} suggests that the discrete Hamilton--Jacobi equation is a generalization of the Riccati equation.
  We are interested in exploring possible further links implied by the generalization.
  \smallskip
\item {\em Galerkin discrete optimal control problems}:
  The Galerkin discrete control Hamiltonians may be considered to be a means of formulating discrete optimal control problems with higher-order of approximation to a continuous-time optimal control problem.
  This idea generalizes the Runge--Kutta discretizations of optimal control problems (see, e.g., \citet{Ha2000} and references therein).
  In fact, \citet{LeZh2011} showed that their method recovers the SPRK (symplectic-partitioned Runge--Kutta) method.
  Therefore, this approach is expected to provide structure-preserving higher-order numerical methods for optimal control problems.
\end{itemize}

\section*{Acknowledgments}
This work was partially supported by NSF grants DMS-604307, DMS-0726263, DMS-0907949, and DMS-1010687.
We would like to thank the referees, Jerrold Marsden, Harris McClamroch, Matthew West, Dmitry Zenkov, and Jingjing Zhang for helpful discussions and comments.

\appendix

\section{Discrete Linear Hamiltonian Systems}
\label{asec:LinearDiscreteHamiltonianSystems}
\subsection{Discrete Linear Hamiltonian Systems}
\label{assec:LinearDiscreteHamiltonianSystems}
Suppose that the configuration space $Q$ is an $n$-dimensional vector space, and that the discrete Hamiltonian $H_{\rm d}^{+}$ or $H_{\rm d}^{-}$ is quadratic as in Eq.~\eqref{eq:QDH}.
Also assume that the corresponding discrete Hamiltonian map $\tilde{F}_{L_{\rm d}}: (q_{k}, p_{k}) \mapsto (q_{k+1}, p_{k+1})$ is invertible.
Then, the discrete Hamilton's equations~\eqref{eq:RDHEq} or \eqref{eq:LDHEq} reduce to the discrete linear Hamiltonian system
\begin{equation}
  \label{eq:DLS}
  z_{k+1} = A_{L_{\rm d}} z_{k},
\end{equation}
where $z_{k} \in \R^{2n}$ is a coordinate expression for $(q_{k}, p_{k}) \in Q \oplus Q^{*}$ and $A_{L_{\rm d}}: Q \oplus Q^{*} \to Q \oplus Q^{*}$ is the matrix representation of the map $\tilde{F}_{L_{\rm d}}$ under the same basis.
Since $\tilde{F}_{L_{\rm d}}$ is symplectic, $A_{L_{\rm d}}$ is an $2n \times 2n$ symplectic matrix, i.e.,
\begin{equation}
  A_{L_{\rm d}}^{T} \mathbb{J} A_{L_{\rm d}} = \mathbb{J},
\end{equation}
where the matrix $\mathbb{J}$ is defined by
\begin{equation*}
  \mathbb{J} \defeq
  \begin{pmatrix}
    0 & I \\
    -I & 0
  \end{pmatrix}
\end{equation*}
with $I$ the $n \times n$ identity matrix.

\subsection{Lagrangian Subspaces and Lagrangian Affine Spaces}
\label{ssec:LagrangianSubAffinceSpaces}
Let us recall the definition of a Lagrangian subspace:
\begin{definition}
  Let $V$ be a symplectic vector space with the symplectic form $\Omega$.
  A subspace $\mathcal{L}$ of $V$ is said to be {\em Lagrangian} if $\Omega(v,w) = 0$ for any $v, w \in \mathcal{L}$ and $\dim\mathcal{L} = \dim V/2$.
\end{definition}
We introduce the following definition for later convenience:
\begin{definition}
  \label{def:LagrangianAffineSpace}
  A subset $\tilde{\mathcal{L}}_{b}$ of a symplectic vector space $V$ is called a {\em Lagrangian affine space} if $\tilde{\mathcal{L}}_{b} = b + \mathcal{L}$ for some element $b \in V$ and a Lagrangian subspace $\mathcal{L} \subset V$.
\end{definition}
The following fact is well-known~(see, e.g., \citet[Theorem~6 on p.~417]{Ju1997}):
\begin{proposition}
  Let $\mathcal{L}$ be a Lagrangian subspace of $V$ and $\Phi: V \to V$ be a symplectic transformation.
  Then, for any $k \in \N$, the image of $\mathcal{L}$ under the $k$-fold composition of $\Phi$, i.e.,
  \begin{equation*}
    \Phi^{k}(\mathcal{L}) = \underbrace{ \Phi \circ \dots \circ \Phi }_{k}(\mathcal{L})
  \end{equation*}
  is also a Lagrangian subspace of $V$.
\end{proposition}
A similar result holds for Lagrangian affine spaces:
\begin{proposition}
  \label{prop:InvLagAffineSp}
  Let $\tilde{\mathcal{L}}_{b} = b + \mathcal{L}$ be a Lagrangian affine space of $V$ and $\Phi: V \to V$ be a symplectic transformation.
  Then $\Phi^{k}( \tilde{\mathcal{L}}_{b} )$ is also a Lagrangian affine space of $V$ for any $k \in \N$.
  More explicitly, we have
  \begin{equation*}
    \Phi^{k}( \tilde{\mathcal{L}}_{b} ) = \Phi^{k}b + \Phi^{k}(\mathcal{L}).
  \end{equation*}
\end{proposition}
\begin{proof}
  Follows from a straightforward calculation.
\end{proof}

\subsection{Generating Functions}
Now, consider the case where $V = Q \oplus Q^{*}$ to apply the results from Section~\ref{ssec:LagrangianSubAffinceSpaces} to the setting in Section~\ref{assec:LinearDiscreteHamiltonianSystems}.
This is a symplectic vector space with the symplectic form $\Omega: (Q \oplus Q^{*}) \times (Q \oplus Q^{*}) \to \R$ defined by
\begin{equation*}
  \Omega: (v, w) \mapsto v^{T} \mathbb{J} w.
\end{equation*}
The key result here regarding Lagrangian subspaces on $Q \oplus Q^{*}$ is the following:
\begin{proposition}
  \label{prop:mathcalL-dS}
  A Lagrangian subspace of $Q \oplus Q^{*}$ that is transversal to $\{0\} \oplus Q^{*}$ is the graph of an exact one-form on $Q$, i.e., $\mathcal{L} = \Graph dS$ for some function $S: Q \to \R$ which has the form
  \begin{equation}
    \label{eq:S(q)}
    S(q) = \frac{1}{2}\ip{A q}{q} + C
  \end{equation}
  with some symmetric linear map $A: Q \to Q^{*}$ and an arbitrary real scalar constant $C$.
  Moreover, the correspondence between the Lagrangian subspaces and such functions (modulo the constant term) is one-to-one.
\end{proposition}
\begin{proof}
  First, recall that a Lagrangian submanifold of $T^{*}Q$ that projects diffeomorphically onto $Q$ is the graph of a closed one-forms on $Q$~(see \citet[Proposition~5.3.15 and the subsequent paragraph on p.~410]{AbMa1978}).
  In our case, $Q$ is a vector space, and so the cotangent bundle $T^{*}Q$ is identified with the direct sum $Q \oplus Q^{*}$.
  Now, a Lagrangian subspace of $Q \oplus Q^{*}$ that is transversal to $\{0\} \oplus Q^{*}$ projects diffeomorphically onto $Q$, and so is the graph of a closed one-form.
  Then, by the Poincar\'e lemma, it follows that any such Lagrangian subspace $\mathcal{L}$ is identified with the graph of an exact one-form $dS$ with some function $S$ on $Q$, i.e., $\mathcal{L} = \Graph dS$.

  However, as shown in, e.g., \citet[][Theorem~3 on p.~233]{Ju1997}, the space of Lagrangian subspaces that are transversal to $\{0\} \oplus Q^{*}$ is in one-to-one correspondence with the space of all symmetric maps $A: Q \to Q^{*}$, with the correspondence given by $\mathcal{L} = \Graph A$.
  Hence, $\Graph dS = \Graph A$, or more specifically, 
  \begin{equation*}
    dS(q) = A_{ij}q^{j}\,dq^{i}.
  \end{equation*}
  This implies that $S$ has the form
  \begin{equation*}
    S(q) = \frac{1}{2}A_{ij} q^{i} q^{j} + C,
  \end{equation*}
  with an arbitrary real scalar constant $C$.
\end{proof}

\begin{corollary}
  \label{cor:GenFcn-LagAffineSp}
  Let $\tilde{\mathcal{L}}_{z_{0}} = z_{0} + \mathcal{L}$ be a Lagrangian affine space, where $z_{0} = (q_{0}, p_{0})$ is an element in $Q \oplus Q^{*}$ and $\mathcal{L}$ is a Lagrangian subspace of $Q \oplus Q^{*}$ that is transversal to $\{0\} \oplus Q^{*}$.
  Then, $\tilde{\mathcal{L}}_{z_{0}}$ is the graph of an exact one-form $d\tilde{S}$ with a function $\tilde{S}: Q \to \R$ of the form
  \begin{equation}
    \label{eq:tildeS}
    \tilde{S}(q) = \frac{1}{2}\ip{A q}{q} + \ip{p_{0} - A q_{0}}{q} + C,
  \end{equation}
  with some symmetric linear map $A: Q \to Q^{*}$ and an arbitrary real scalar constant $C$.
\end{corollary}
\begin{proof}
  From Proposition~\ref{prop:mathcalL-dS}, there exists a function $S: Q \to \R$ of the form Eq.~\eqref{eq:S(q)} with some symmetric linear map $A: Q \to Q^{*}$ such that $\mathcal{L} = \Graph dS = \Graph A$.
  Let us define $\tilde{S}: Q \to R$ by
  \begin{equation*}
    \tilde{S}(q) \defeq S(q - q_{0}) + \ip{p_{0}}{q},
  \end{equation*}
  from which Eq.~\eqref{eq:tildeS} follows by direct calculation.
  Then,
  \begin{equation*}
    d\tilde{S}(q) = A(q - q_{0}) + p_{0}.
  \end{equation*}
  and thus
  \begin{align*}
    \Graph d\tilde{S}
    &= \{ (q, d\tilde{S}(q) )\ |\ q \in Q \}
    \\
    &= \setdef{ \parentheses{q, A(q - q_{0}) + p_{0}}}{ q \in Q }
    \\
    &= (q_{0}, p_{0}) + \setdef{ \parentheses{q - q_{0}, A(q - q_{0})} }{ q \in Q }
    \\
    &= z_{0} + \Graph A
    \\
    &= z_{0} + \mathcal{L}
    \\
    &= \tilde{\mathcal{L}}_{z_{0}}. \qedhere
  \end{align*}
\end{proof}

\bibliography{DiscreteHJ}

\begin{thebibliography}{34}
\providecommand{\natexlab}[1]{#1}
\providecommand{\url}[1]{\texttt{#1}}
\expandafter\ifx\csname urlstyle\endcsname\relax
  \providecommand{\doi}[1]{doi: #1}\else
  \providecommand{\doi}{doi: \begingroup \urlstyle{rm}\Url}\fi

\bibitem[Abraham and Marsden(1978)]{AbMa1978}
R.~Abraham and J.~E. Marsden.
\newblock \emph{Foundations of Mechanics}.
\newblock Addison--Wesley, 2nd edition, 1978.

\bibitem[Ammar and Martin(1986)]{AmMa1986}
G.~Ammar and C.~Martin.
\newblock The geometry of matrix eigenvalue methods.
\newblock \emph{Acta Applicandae Mathematicae}, 5\penalty0 (3):\penalty0
  239--278, 1986.

\bibitem[Arnold(1989)]{Ar1991}
V.~I. Arnold.
\newblock \emph{Mathematical Methods of Classical Mechanics}.
\newblock Springer, 1989.

\bibitem[Bellman(1971)]{Be1971}
R.~Bellman.
\newblock \emph{Introduction to the Mathematical Theory of Control Processes},
  volume~2.
\newblock Academic Press, 1971.

\bibitem[Bellman(1972)]{Be1972}
R.~Bellman.
\newblock \emph{Dynamic programming}.
\newblock Princeton University Press, 1972.

\bibitem[Bertsekas(2005)]{Be2005}
D.~P. Bertsekas.
\newblock \emph{Dynamic Programming and Optimal Control}, volume~1.
\newblock Athena Scientific, 2005.

\bibitem[Bloch(2003)]{Bl2003}
A.~M. Bloch.
\newblock \emph{Nonholonomic Mechanics and Control}.
\newblock Springer, 2003.

\bibitem[Brockett(1981)]{Br1981}
R.~W. Brockett.
\newblock Control theory and singular {R}iemannian geometry.
\newblock In P.~J. Hilton and G.~S. Young, editors, \emph{New Directions in
  Applied Mathematics}, pages 11--27. Springer, 1981.

\bibitem[Cadzow(1970)]{Ca1970}
J.~A. Cadzow.
\newblock Discrete calculus of variations.
\newblock \emph{International Journal of Control}, 11\penalty0 (3):\penalty0
  393--407, 1970.

\bibitem[Cari\~nena et~al.(2010)Cari\~nena, Gracia, Marmo, Mart\'inez, Mun\~oz
  Lecanda, and Rom\'an-Roy]{CaGrMaMaMuRo2010}
J.~F. Cari\~nena, X.~Gracia, G.~Marmo, E.~Mart\'inez, M.~C. Mun\~oz Lecanda,
  and N.~Rom\'an-Roy.
\newblock Geometric {H}amilton--{J}acobi theory for nonholonomic dynamical
  systems.
\newblock \emph{International Journal of Geometric Methods in Modern Physics},
  7\penalty0 (3):\penalty0 431--454, 2010.

\bibitem[Channell and Scovel(1990)]{ChSc1990}
P.~J. Channell and C.~Scovel.
\newblock Symplectic integration of {H}amiltonian systems.
\newblock \emph{Nonlinearity}, 3\penalty0 (2):\penalty0 231--259, 1990.

\bibitem[de~Le\'on et~al.(2010)de~Le\'on, Marrero, and Mart\'in~de
  Diego]{LeMaMa2010}
M.~de~Le\'on, J.~C. Marrero, and D.~Mart\'in~de Diego.
\newblock Linear almost {P}oisson structures and {H}amilton--{J}acobi equation.
  {A}pplications to nonholonomic mechanics.
\newblock \emph{Journal of Geometric Mechanics}, 2\penalty0 (2):\penalty0
  159--198, 2010.

\bibitem[Elnatanov and Schiff(1996)]{ElSc1996}
N.~A. Elnatanov and J.~Schiff.
\newblock The {H}amilton--{J}acobi difference equation.
\newblock \emph{Functional Differential Equations}, 3\penalty0 (279--286),
  1996.

\bibitem[Ge and Marsden(1988)]{GeMa1988}
Z.~Ge and J.~E. Marsden.
\newblock Lie--{P}oisson {H}amilton--{J}acobi theory and {L}ie--{P}oisson
  integrators.
\newblock \emph{Physics Letters A}, 133\penalty0 (3):\penalty0 134--139, 1988.

\bibitem[Gelfand and Fomin(2000)]{GeFo2000}
I.~M. Gelfand and S.~V. Fomin.
\newblock \emph{Calculus of Variations}.
\newblock Dover, 2000.

\bibitem[Goldstein et~al.(2001)Goldstein, Poole, and Safko]{GoPoSa2001}
H.~Goldstein, C.~P. Poole, and J.~L. Safko.
\newblock \emph{Classical Mechanics}.
\newblock Addison Wesley, 3rd edition, 2001.

\bibitem[Guibout and Bloch(2004)]{GuBl2004}
V.~Guibout and A.~M. Bloch.
\newblock A discrete maximum principle for solving optimal control problems.
\newblock In \emph{43rd IEEE Conference on Decision and Control}, volume~2,
  pages 1806--1811 Vol.2, 2004.

\bibitem[Hager(2000)]{Ha2000}
W.~W. Hager.
\newblock Runge--{K}utta methods in optimal control and the transformed adjoint
  system.
\newblock \emph{Numerische Mathematik}, 87\penalty0 (2):\penalty0 247--282,
  2000.

\bibitem[Iglesias-Ponte et~al.(2008)Iglesias-Ponte, de~Le\'on, and Mart\'in~de
  Diego]{IgLeMa2008}
D.~Iglesias-Ponte, M.~de~Le\'on, and D.~Mart\'in~de Diego.
\newblock Towards a {H}amilton--{J}acobi theory for nonholonomic mechanical
  systems.
\newblock \emph{Journal of Physics A: Mathematical and Theoretical},
  41\penalty0 (1), 2008.

\bibitem[Jalnapurkar et~al.(2006)Jalnapurkar, Leok, Marsden, and
  West]{JaLeMaWe2006}
S.~M. Jalnapurkar, M.~Leok, J.~E. Marsden, and M.~West.
\newblock Discrete {R}outh reduction.
\newblock \emph{Journal of Physics A: Mathematical and General}, 39\penalty0
  (19):\penalty0 5521--5544, 2006.

\bibitem[Jordan and Polak(1964)]{JoPo1964}
B.~W. Jordan and E.~Polak.
\newblock Theory of a class of discrete optimal control systems.
\newblock \emph{Journal of Electronics and Control}, 17:\penalty0 694--711,
  1964.

\bibitem[Jurdjevic(1997)]{Ju1997}
V.~Jurdjevic.
\newblock \emph{Geometric control theory}.
\newblock Cambridge University Press, Cambridge, 1997.

\bibitem[Lall and West(2006)]{LaWe2006}
S.~Lall and M.~West.
\newblock Discrete variational {H}amiltonian mechanics.
\newblock \emph{Journal of Physics A: Mathematical and General}, 39\penalty0
  (19):\penalty0 5509--5519, 2006.

\bibitem[Lanczos(1986)]{La1986}
C.~Lanczos.
\newblock \emph{The Variational Principles of Mechanics}.
\newblock Dover, 4th edition, 1986.

\bibitem[Leok(2004)]{Le2004}
M.~Leok.
\newblock \emph{Foundations of Computational Geometric Mechanics}.
\newblock PhD thesis, California Institute of Technology, 2004.

\bibitem[Leok and Zhang(2011)]{LeZh2011}
M.~Leok and J.~Zhang.
\newblock Discrete {H}amiltonian variational integrators.
\newblock \emph{IMA Journal of Numerical Analysis}, 2011.

\bibitem[Leok et~al.()Leok, Ohsawa, and Sosa]{DiracHJ}
M.~Leok, T.~Ohsawa, and D.~Sosa.
\newblock {H}amilton--{J}acobi theory for degenerate {L}agrangian systems with
  constraints.
\newblock \emph{in preparation}.

\bibitem[Leok et~al.(2004)Leok, Marsden, and Weinstein]{LeMaWe2004}
M.~Leok, J.~E. Marsden, and A.~Weinstein.
\newblock A discrete theory of connections on principal bundles.
\newblock \emph{Preprint}, 2004.

\bibitem[Marsden and Ratiu(1999)]{MaRa1999}
J.~E. Marsden and T.~S. Ratiu.
\newblock \emph{Introduction to Mechanics and Symmetry}.
\newblock Springer, 1999.

\bibitem[Marsden and West(2001)]{MaWe2001}
J.~E. Marsden and M.~West.
\newblock Discrete mechanics and variational integrators.
\newblock \emph{Acta Numerica}, pages 357--514, 2001.

\bibitem[Moser and Veselov(1991)]{MoVe1991}
J.~Moser and A.~P. Veselov.
\newblock Discrete versions of some classical integrable systems and
  factorization of matrix polynomials.
\newblock \emph{Communications in Mathematical Physics}, 139\penalty0
  (2):\penalty0 217--243, 1991.

\bibitem[Ohsawa and Bloch(2009)]{OhBl2009}
T.~Ohsawa and A.~M. Bloch.
\newblock Nonholonomic {H}amilton--{J}acobi equation and integrability.
\newblock \emph{Journal of Geometric Mechanics}, 1\penalty0 (4):\penalty0
  461--481, 2009.

\bibitem[Suris(2003)]{Su2003}
Y.~B. Suris.
\newblock \emph{The problem of integrable discretization: {H}amiltonian
  approach}.
\newblock Birkh\"auser, Basel, 2003.

\bibitem[Suris(2004)]{Su2004}
Y.~B. Suris.
\newblock Discrete {L}agrangian models.
\newblock In \emph{Discrete Integrable Systems}, volume 644 of \emph{Lecture
  Notes in Physics}, pages 111--184. Springer, 2004.

\end{thebibliography}
\bibliographystyle{plainnat}

\end{document}